\documentclass[reqno,a4paper]{amsart}
\usepackage{amsmath,amsthm,amssymb,tabularx}

\theoremstyle{plain}
\newtheorem{tm}{Theorem}[section]
\newtheorem{lm}[tm]{Lemma}
\newtheorem{cor}[tm]{Corollary}
\newtheorem{prop}[tm]{Proposition}

\theoremstyle{definition}
\newtheorem{definition}[tm]{Definition}
\newtheorem{notations}[tm]{Notations}
\newtheorem{ex}[tm]{Example}
\newtheorem{example}[tm]{Example}
\newtheorem{rem}[tm]{Remark}
\newtheorem{nots}[tm]{Notations}
\allowdisplaybreaks

\newcommand{\beq}{\begin{equation}}
\newcommand{\eeq}{\end{equation}}
\newcommand{\bga}{\begin{gather*}}
\newcommand{\ega}{\end{gather*}}
\newcommand{\bal}{\begin{align*}}
\newcommand{\eal}{\end{align*}}
\newcommand{\bit}{\begin{itemize}}
\newcommand{\eit}{\end{itemize}}
\newcommand{\btm}{\begin{tm}}
\newcommand{\etm}{\end{tm}}
\newcommand{\blm}{\begin{lm}}
\newcommand{\elm}{\end{lm}}
\newcommand{\bcor}{\begin{cor}}
\newcommand{\ecor}{\end{cor}}
\newcommand{\bprop}{\begin{prop}}
\newcommand{\eprop}{\end{prop}}
\newcommand{\bex}{\begin{ex}}
\newcommand{\eex}{\end{ex}}
\newcommand{\bpr}{\begin{proof}}
\newcommand{\epr}{\end{proof}}
\newcommand{\brem}{\begin{rem}}
\newcommand{\erem}{\end{rem}}
\newcommand{\bdf}{\begin{definition}}
\newcommand{\edf}{\end{definition}}
\newcommand{\bnots}{\begin{nots}}
\newcommand{\enots}{\end{nots}}

\def\C{\mathbb{C}}
\def\N{\mathbb{N}}
\def\R{\mathbb{R}}

\def\Id{\mathbb I}

\def\e{\varepsilon}
\let\a\alpha

\let\cal\mathcal

\def\id{{\rm id}}

\def \le {\leqslant}
\def \ge {\geqslant}
\let \<\langle
\let \>\rangle
\let\phi\varphi
\let\kappa\varkappa

\def\cst{$C^*$}
\def\Mhat{\widehat M}
\let\hat\widehat

\def\QSI{\cal {QSI}}
\def\G{\mathbb{G}}

\begin{document}
\author{Yulia Kuznetsova}
\title{Duals of quantum semigroups with involution}
\thanks{This work was partially supported by the Simons Foundation grant 346300 and the Polish Government MNiSW 2015-2019 matching fund. The author was also supported by the travel grant PHC Star 2016 36618SE of the French Ministry of Foreign Affairs.}
\address{University of Bourgogne Franche-Comt\'e, 16 route de Gray, 25030 Besan\c con, France}
\email{yulia.kuznetsova@univ-fcomte.fr}

\begin{abstract}
We define a category $\QSI$ of quantum semigroups with involution which carries a corep\-re\-sen\-ta\-tion-based duality map $M\mapsto \widehat M$. Objects in $\QSI$ are von Neumann algebras with comultiplication and coinvolution, we do not suppose the existence of a Haar weight or of a distinguished spatial realisation.
In the case of a locally compact quantum group $\G$, the duality $\;\widehat{\ }\;$ in $\QSI$ recovers
the universal duality of Kustermans: $\widehat{L^\infty(\G)} = C_0^u(\hat \G)^{**}= \widehat{ C_0^u(\G)^{**}}$, and $\widehat{L^\infty(\hat\G)} = C_0^u(\G)^{**} = \widehat{ C_0^u(\hat\G)^{**}}$. Other various examples are given.
\end{abstract}

\maketitle

\setcounter{tocdepth}{1}

\section{Introduction}

In the theory of locally compact groups, a crucial role is played by the Haar measure; in the theory of locally compact quantum groups, its place it taken by a
pair of left and right Haar weights. It is known that in the classical case, the existence of a Haar measure is a theorem, whereas in the quantum case this is a part of axioms.

The dual of a locally compact quantum group in the setting of Kustermans and Vaes \cite{kust-vaes} is defined via an explicit construction involving the Haar weight; in the setting of Woronowicz \cite{woron}, the spatial realisation of the multiplicative unitary gives readily the corresponding pair of quantum group algebras in duality.

This should be compared with the classical Pontryagin duality: we consider first the group $\widehat G$ of continuous characters of $G$, which is meaningful for any topological group, and then prove that if $G$ is locally compact abelian, then so is $\widehat G$, and $\widehat{\widehat G}$ is isomorphic to $G$.

It is known that in the quantum setting, one can define a dual in a measure-independent way, but arriving then to a different (bigger) object called the {\it universal dual.} This was done by Kirchberg \cite{kirchberg} in the particular case of Kac algebras, extended to locally compact quantum groups (LCQG) by Kustermans \cite{kustermans}, and later realized by So\l tan and Woronowicz \cite{sol-wor} in the setting of multiplicative unitaries.

The problem is however that this approach does not give any existence theorem except for the case of LCQG or multiplicative unitaries respectively. 
In general, there is no guarantee that the dual of the universal dual
exists, and thus we cannot speak of a Pontryagin-like isomorphism between an algebra and its bidual. 

The present paper introduces a category $\QSI$ of {\it quantum semigroups with involution} which solves this problem. For every $M\in \QSI$, we define its dual object $\Mhat\in\QSI$, so that it extends the universal duality of Kustermans: if $L^\infty(\G)$ is a von Neumann algebraic LCQG and $\;\widehat{\ }\;$ denotes the dual in $\QSI$, then $\widehat{L^\infty(\G)} = C_0^u(\hat \G)^{**}= \widehat{ C_0^u(\G)^{**}}$, and $\widehat{L^\infty(\hat\G)} = C_0^u(\G)^{**} = \widehat{ C_0^u(\hat\G)^{**}}$. Moreover, the construction produces first the $C^*$-algebras $C^u_0(\hat \G)$ and $C^u_0(\G)$ and only then we pass to their universal enveloping von Neumann algebras (identified with the linear biduals).

An object in $\QSI$ is a von Neumann algebra $M$ with comultiplication $\Delta: M\to M\otimes M$ and a densely defined coinvolution $S:D(S)\subset M\to M$, which satisfy certain natural axioms (see Section~\ref{sec-defs}). This class includes all LCQG, but a more general example would be the algebra $C(P)$ of continuous functions on a compact semitopological semigrouop with involution $P$ (Example \ref{ex-semigroup}).

The structure on $M$ induces on the predual $M_*$ a structure of a Banach algebra, with an involution defined on a subalgebra $M_{**}$. By passing to a quotient of $M$, we can guarantee that $M_{**}$ is dense in $M_*$. We consider representations of $M_*$ which are involutive on $M_{**}$, and, as usual, call such a representation unitary if it is generated by a unitary $U\in M\otimes B(H)$ (Section \ref{sec-ideal}).

The central point is to consider the ideal $M_{**}^0\subset M_*$ defined as the common kernel of all {\it irreducible non-unitary\/} representations. It appears that this ideal contains all necessary information on the dual algebra.

We set $\widehat M = C^*(M_{**}^0)^{**}$. By definition (except for degenerate cases when we set $\Mhat=\{0\}$), irreducible representations of this algebra corrrespond bijectively to the unitary irreducible representations of $M_*$; by disintegration, this bijection extends also to reducible repesentations, see Section~\ref{sec-ideal}.

Every dual algebra carries a structure of a quantum semigroup with involution. Moreover, it is a Hopf algebra in a certain sense (Subsection \ref{sec-antipode}), what suggest to consider it as a quantum {\it group} rather than a semigroup.

If $M=\widehat N$ is a dual of some $N\in \QSI$ and is either commutative or cocommutative, then it is isomorphic to $C_0(G)^{**}$ or $C^*(G)^{**}$ respectively for a classical locally compact group $G$. Moreover, in all examples known (Sections \ref{sec-LCQG} and \ref{sec-examples}) the dual algebra $\Mhat$ coincides with a dual of a LCQG, so in particular $\Mhat$ is isomorphic to the third dual of $M$. This allows to conjecture that $\Mhat\simeq\widehat{\widehat{\Mhat}}$ in general.

This article is a development of \cite{haar}, where the coinvolution of $M$ (termed antipode there) was supposed to be bounded.

\section{Definitions and notations}\label{sec-defs}

Our main objects are von Neumann bialgebras, which we allow to be just zero spaces. Below, $\otimes$ denotes the von Neumann algebraic tensor product, and $\Id$ the identity operator.

\begin{definition}
A {\it von Neumann bialgebra} $M$ is a von Neumann algebra endowed with a comultiplication $\Delta: M\to M\otimes M$ which is a normal coassociative *-homomorphism: $\Delta*=(*\otimes*)\Delta$ and $(\Id\otimes\Delta)\Delta = (\Delta\otimes\Id)\Delta$.
\end{definition}

Note that we do not require $\Delta$ to be unital.

It is known that in the assumptions above the predual $M_*$ of $M$ is a completely contractive Banach algebra.

In addition to this structure, we postulate the existence of a {\it coinvolution} $S$ on $M$. As it will be seen, the requirements on $S$ are weaker than on an atipode, and this notion corresponds rather to an {\it involution} of a semigroup than to a group inverse. For this reason the author prefers to use the term coinvolution and not antipode as it was done in \cite{haar}.

\begin{definition}\label{def-coinvolution}
Let $M$ be a von Neumann bialgebra. A linear map $S: D(S)\subset M\to M$ is called a {\it proper coinvolution} if it satisfies the following conditions:
\begin{enumerate}
\item $D(S)$ is $\sigma$-weakly dense in $M$;
\item $D(S)$ is closed under multiplication and $S: D(S)\to M$ is an anti-homomorphism;
\item $(*S)(D(S))\subset D(S)$ and $(*S)^2=\Id_{D(S)}$;
\item if $\mu,\nu\in M_*$ are such that $\mu\circ S$ and $\nu\circ S$ extend to normal functionals on $M$, then for all $x\in D(S)$ holds $\Delta (x)(\nu\circ S\otimes\mu\circ S) = \Delta S(x)(\mu\otimes\nu)$.
\end{enumerate}
A von Neumann bialgebra equipped with a proper coinvolution is called a {\it quantum semigroup with involution}.
\end{definition}

The last formula is a replacement of the identity $\theta(S\otimes S)\Delta=\Delta S$ which may have no sense a priori. The antipode of a locally compact quantum group satisfies these conditions, see \cite[Lemma 5.25]{kust-vaes}. Note that we do not require $S$ to be closeable.

\begin{example}\label{ex-semigroup}
Let $P$ be a compact semitopological semigroup with involution. Recall that an involution on a semigroup is a map $*:P\to P$ such that $(x^*)^*=x$ and $(xy)^*=y^*x^*$ for all $x,y\in P$. Let $SC(P\times P)$ be the space of separately continuous functions on $P\times P$, then we get a natural map $\check \Delta: C(P)\to SC(P\times P)$, $\check\Delta(f)(s,t)=f(st)$. Set $M=C(P)^{**}$, with the usual structure of a von Neumann algebra. Since $SC(P\times P)$ is canonically imbedded into $M\otimes M$, $\check \Delta$ can be viewed as a map from $C(P)$ to $M\otimes M$, and it is known \cite{das-mroz} that it is a unital *-homomorphism and as such can be extended by normality to $M$, so that the extension $\Delta$ satisfies $(\Delta\otimes \id)\check\Delta = (\id\otimes\Delta)\check\Delta$. Altogether, this implies that $\Delta$ is a comultiplication on $M$. Set $D(S)=C(P)$ and $(Sf)(t)=f(t^*)$ for $f\in C(P)$, $t\in P$. It is easily seen that $S$ satisfied conditions (1)--(4), so that $M=C(P)^{**}$ is a quantum semigroup with involution.
\end{example}

\begin{definition}\label{def-involution}
For every $\mu\in M_*$, define $\bar\mu\in M_*$ by $\bar\mu(a)=\overline{\mu(a^*)}$, $a\in M$. Let $M_{**}$ be the subspace of all $\mu\in M_*$ such that $\bar\mu\circ S$ extends to a bounded normal functional on $M$. We will denote by $\mu^*$ this extension, so that $\mu^*(x) = \bar\mu(Sx)$ for $x\in D(S)$.
\end{definition}

Note that in $M_*$ as a predual space, $M_{**}$ is $\sigma(M_*,M)$-dense iff it is norm dense.

\bprop
If $M_{**}$ is dense in $M_*$ then $S$ is closeable.
\eprop
\bpr
Suppose that $M_{**}$ is dense in $M_*$, that is, separates the points of $M$. Towards a contradiction, suppose also that $S$ is not closeable. Then
there exist $a_n\to a$, $b_n\to a$ such that $\a=\lim Sa_n\ne \lim Sb_n=\beta$. Then for every $\mu\in M_{**}$ we have:
$\bar\mu(\a) = \lim \bar\mu(Sa_n) = \lim\mu^*(a_n) = \mu^*(a) = \lim\mu^*(b_n) = \lim \bar\mu(Sb_n) = \bar\mu(\beta)$. Thus $\mu(\a^*-\beta^*)=0$ for all $\mu\in M_{**}$ while $\a-\beta\ne0$, and we arrive at a contradiction.
\epr

\bprop
On $M_*$, $\mu\mapsto \bar\mu$ is a homomorphism. If $S^2=\Id$, it is involutive.
\eprop
\bpr
$$
(\bar\mu\cdot\bar\nu)(x) = (\bar\mu\otimes\bar\nu)(\Delta(x)) = \overline{(\mu\otimes\nu)(\Delta(x)^*)}
 = \overline{(\mu\otimes\nu)(\Delta(x^*))} = \overline{(\mu\cdot\nu)(x^*)} = (\overline{\mu\cdot\nu})(x).
$$
\epr

\bprop
$M_{**}$ is a subalgebra in $M_*$, and $(\mu\nu)^* = \nu^*\mu^*$ for all $\mu,\nu\in M_{**}$.
With the norm $\|\mu\|_* = \max\big( \|\mu\|, \|\mu^*\| \big)$, $M_{**}$ is a Banach *-algebra.
\eprop
\bpr
For $\mu,\nu\in M_{**}$ and $x\in D(S)$, we have with \ref{def-coinvolution}(4):
\begin{align*}
\overline{\mu\nu}\circ S(x) &=
 (\bar\mu\cdot\bar\nu)(Sx) = (\bar\mu\otimes\bar\nu) \big(\Delta(Sx)\big)
 = (\bar\nu\circ S\otimes\bar\mu\circ S) \big(\Delta(x)\big)
 = (\nu^*\mu^*) (x).
\end{align*}
The first statement follows.

It is immediate that $\|\mu\nu\|_*\le \|\mu_*\|\nu\|_*$. It remains to prove that $(M_{**},\|\cdot\|_*)$ is complete. If $(\mu_n)\subset M_{**}$ is a Cauchy sequence with respect to $\|\cdot\|_*$, then in particular $\mu_n\to\mu$, $\mu_n^*\to\nu$ with respect to $\|\cdot\|$ for some $\mu,\nu\in M_*$.
For every $a\in D(S)$, we have $\bar\mu(Sa) = \overline{\mu\big((Sa)^*\big)} = \lim \overline{\mu_n\big((Sa)^*\big)}
 = \lim \mu_n^*(a) = \nu(a)$, thus $\mu\in M_{**}$ and $\mu^*=\nu$, moreover $\mu_n\to\mu$ in $(M_{**},\|\cdot\|_*)$.
\epr

Recall that on $M_*$, there is a two-sided action of $M$. For $\mu\in M_*$, $a\in M$ one sets $(\mu.a)(b) = \mu(ab)$ and $(a.\mu)(b)=\mu(ba)$, $b\in M$.
As in the case of a locally compact group \cite[Lemma 2.1]{kustermans}, $M_{**}$ is stable under multiplication by $D(S)$:

\blm\label{perp-module-MS}
$M_{**}$ is a two-sided module over $D(S)$.
\elm
\bpr
Let $\mu$ be in $M_{**}$, $a\in D(S)$. Then for $b\in D(S)$ we have
\begin{align*}
(\overline{a.\mu\strut})(Sb) &= \overline{(a.\mu)( (Sb)^*)}
 = \overline{\mu( (Sb)^* a)}  = \overline{\mu( (Sb)^* ( S(Sa)^*)^*)}
\\& = \overline{\mu( (S(b(Sa)^*)^*)} = \mu^*( b(Sa)^* ) = ( (Sa)^*. \mu^* )( b ).
\end{align*}
This means that $\overline{(a.\mu)}\circ S$ extends to the bounded functional $(Sa)^*. \mu^*=(a.\mu)^*$ on $M$, so $a.\mu\in M_{**}$. The proof for the right action is similar.
\epr

\begin{definition}\label{def-R}
Let $M_{**}$ be dense in $M_*$. Set
$$
D(\widetilde S) = \{ a\in M: \exists C>0 \text{ such that } |\mu^*(a)|\le C\|\mu\| \text{ for all }\mu\in M_{**}\}.
$$
For every $a\in D(\widetilde S)$, define $\widetilde S(a)\in M$ by $\mu\big(\widetilde S(a)\big) = \overline{\mu^*(a)}$.
\end{definition}

\bprop
Let $M_{**}$ be dense in $M_*$. Then $\widetilde S$ is a closed operator.
\eprop
\bpr
If $a_n\to a$, $\widetilde S(a_n)\to b$ in $M$, then for every $\mu\in M_{**}$
$\mu^*(a_n) \to \mu^*(a)$; $\mu^*(a_n) = \overline{\mu\big( \widetilde S(a_n)\big) } \to \overline{ \mu(b)}$ as $n\to\infty$. It follows that $\mu^*(a) = \overline{ \mu(b)}$ whence $a\in D(\widetilde S)$ and $\widetilde S(a) = b$.
\epr

\brem
Let $M_{**}$ be dense in $M_*$. Then $D(S)\subset D(\widetilde S)$ and $\widetilde S(a)=*S(a)$ for $a\in D(S)$. In particular, $\widetilde S D(S)\subset D(S)$.
\erem

\bprop\label{R-continuous}
Let $M_{**}$ be dense in $M_*$. Then $\widetilde S: D(\widetilde S)\to M$ is $\sigma(M,M_{**})$ continuous.
\eprop
\bpr
Follows from the equality $|\mu(\widetilde S(a))| = |\mu^*(a)|$ for all $\mu\in M_{**}$, $a\in M$.
\epr

\bprop\label{R-hom}
If $M_{**}$ is dense in $M_*$ then:\\
(i) $D(\widetilde S)$ is closed under multiplication and $\widetilde S: D(\widetilde S)\to M_*$ is a homomorphism;\\
(ii) $M_{**}$ is a two-sided module over $D(\widetilde S)$;\\
(iii) $\widetilde S\big( D(\widetilde S)\big) \subset D(\widetilde S)$ and $\widetilde S^2=\Id_{D(\widetilde S)}$.
\eprop
\bpr
(i). First note that $\widetilde S=*S$ is a homomorphism on $D(S)$.

Next we prove that $ab\in D(\widetilde S)$ and $\widetilde S(ab)=\widetilde S(a)\widetilde S(b)$ for $a\in D(\widetilde S)$, $b\in D(S)$. By assumption \ref{def-coinvolution}.(1), there exist $a_i\in D(S)$ such that $a_i\to a$ $\sigma$-weakly. Then $a_ib\to ab$, and by Proposition \ref{R-continuous},
$\widetilde S(a_i)\to \widetilde S(a)$, both $\sigma$-weakly.
For every $\mu\in M_{**}$
\begin{align*}
\overline{\mu^*(ab)} &= \lim \overline{\mu^*(a_ib)} = \lim \mu(\widetilde S(a_ib))| = \lim\mu( \widetilde S(a_i) \widetilde S(b))
 \\&= \lim \widetilde S(b).\mu (\widetilde S(a_i)) = (\widetilde S(b).\mu)(\widetilde S(a)) = \mu(\widetilde S(a)\widetilde S(b)),
\end{align*}
what proves that $ab\in D(\widetilde S)$ and $\widetilde S(ab)=\widetilde S(a)\widetilde S(b)$.

In the same way we prove the equality for $a\in D(S)$, $b\in D(\widetilde S)$.

We can now prove (ii): for $a\in D(\widetilde S)$, $\mu\in M_{**}$, $x\in D(S)$,
\begin{align*}
(\mu. \widetilde S(a))^- ( Sx)
& = \overline{(\mu. \widetilde S(a)) ( (Sx)^* )}
 = \overline{\mu ( \widetilde S(a) (Sx)^* )}
 \\&
 = \overline{\mu ( \widetilde S(a) \widetilde S(x) )}
 = \overline{\mu ( \widetilde S(ax) )}
 = \mu^*(ax)
 = (\mu^*.a)(x).
\end{align*}
This implies that $\mu. \widetilde S(a)\in M_{**}$ and $\big(\mu. \widetilde S(a)\big)^* = \mu^*.a$.
Similarly one shows that $\widetilde S(a).\mu\in M_{**}$.

If now $a,b \in D(\widetilde S)$, then  
for $\mu\in M_{**}$
$$
\mu^*(ab) = (\mu^*.a)(b) = \big(\mu. \widetilde S(a)\big)^*(b) = \overline{\big(\mu .\widetilde S(a)\big)\big(\widetilde S(b)\big)}
 = \overline{\mu\big(\widetilde S(a)\widetilde S(b)\big)};
$$
this shows that $ab\in D(\widetilde S)$ and $\widetilde S(ab) = \widetilde S(a)\widetilde S(b)$.

For (iii), suppose that $a\in D(\widetilde S)$. For $\mu\in M_{**}$,
$$
\mu^*\big( \widetilde S(a)\big) = \overline{ \mu^{**}(a)} = \overline{\mu(a)},
$$
which implies that $\widetilde S(a)\in D(\widetilde S)$ and $\widetilde S(\widetilde S(a)) = a$.
\epr

\bprop
Let $M_{**}$ be dense in $M_*$. Then $\bar S=*\widetilde S$, defined on $D(\widetilde S)$, satisfies the axioms (1)-(4) of a proper coinvolution.
\eprop
\bpr
Clearly $\bar S=S$ on $D(S)$. The axiom (1) holds obviously; (2) and (3) follow immediately from \ref{R-hom}. For (4), suppose that $\mu,\nu\in M_*$ are such that $\mu\circ \bar S$ and $\nu\circ \bar S$ extend to normal functionals on $M$. 
Then $\bar\mu$, $\bar\nu\in M_{**}$, and $\mu\circ \bar S=(\bar\mu)^*$, $\nu\circ \bar S=(\bar\nu)^*$.
For $x\in D(\bar S)=D(\widetilde S)$,
\begin{align*}
(\nu\circ \bar S\otimes\mu\circ \bar S)\big(\Delta (x)\big)
&= \big( (\bar\nu)^*\otimes (\bar\mu)^*\big) \big(\Delta(x)\big) =
 \big( (\bar\nu)^*(\bar\mu)^*\big) (x) = \big( (\bar\mu\bar\nu)^*\big)(x)
 \\&= \overline{ (\overline{\mu\nu})(\widetilde S(x))}
 = (\mu\nu)(\widetilde S(x)^*) = (\mu\otimes\nu)\big(\Delta(\bar S(x))\big).
\end{align*}
\epr

\brem
Let $M_{**}$ be dense in $M_*$. The set $\{\mu\in M_*: \bar\mu\circ \bar S \text{ extends to a normal functional on  }M\}$ is equal to $M_{**}$, and the application of Definition \ref{def-R} to $\bar S$ leads to the same operator $\widetilde S$, now with domain $D(\widetilde S)=D(\bar S)$. Thus, $\bar{\bar S}=\bar S$.
\erem

\bprop\label{S-extended}
Let $M_{**}$ be dense in $M_*$. Then $S$ can be extended to a closed proper coinvolution, with the same space $M_{**}$. From now on, we assume that $S$ is closed and $D(S)=D(\widetilde S)$.
\eprop

\bprop\label{M**perp-ideal}
$M_{**}^\perp$ is a ($\sigma$-weakly closed) ideal in $M$.
\eprop
\bpr
Denote $L=M_{**}^\perp$. Obviously it is closed. By Lemma \ref{perp-module-MS}, if $a\in D(S)$, $b\in L$, $\mu\in M_{**}$, then $\mu(ab)= (a\cdot\mu)(b)=0$, so $D(S)L\subset L$, and we have, with $[\cdot]$ meaning $\sigma$-weak closure in $M$:
$$
ML = [D(S)] L \subset [D(S) L] \subset [L] = L,
$$
and similarly for the other inclusion.
\epr

\begin{notations}\label{def-Z}
Proposition \ref{M**perp-ideal} implies that $M_{**}^\perp = (1-\zeta)M$ for a central projection $\zeta$, which will be denoted also by $\zeta_M$ to indicate the algebra $M$. By general relations in Banach spaces, $[M_{**}] = (M_{**}^\perp)_\perp$, where $[\cdot]$ means the norm closure in $M_*$ and $(\cdot)_\perp$ the annihilator in $M_*$ of a subset of $M$. Denote $M_r=M/M_{**}^\perp$ and let $Q: M\to M_r$ be the quotient map.\\
The isomorphism $Z: M_r\to \zeta M$ is a right inverse for $Q$, and its preadjoint map $Z_*: M_*\to (M_r)_*$ is left inverse to $Q_*:(M_r)_*\to M_*$.
\end{notations}

\bprop\label{quotient-by-perp}
With the quotient structure, $M_r$ is a von Neumann bialgebra, with the predual isomorphic to $[M_{**}]$ (the norm closure of $M_{**}$ in $M_*$).
\eprop
\bpr
In the notations of Remark \ref{def-Z}, $Q(a)=Q(\zeta a)$ for all $a\in M$. We define the comultiplication on $M/M_{**}^\perp$ by $\tilde \Delta(Qa) = (Q\otimes Q)\big(\Delta(a)\big)$, $a\in M$. To show that it is well defined, suppose that $Qa=0$, that is $a\in M_{**}^\perp$. Then for $\mu,\nu\in M_*$
$$
(\mu\otimes\nu)\big((\zeta\otimes \zeta)\Delta(a)\big) = (Z_*\mu\otimes Z_*\nu)\big(\Delta(a)\big)
 = (Z_*\mu\cdot Z_*\nu)(a) =0,
$$
since $Z_*\mu\cdot Z_*\nu\in M_{**}$. It follows that
$$
(Q\otimes Q)\big(\Delta(a)\big) = (Q\otimes Q)(\zeta\otimes \zeta)\big(\Delta(a)\big) = 0.
$$
It is immediate to verify that $\tilde\Delta$ is a comultiplication indeed. Note that by definiton, $Q$ is a coalgebra morphism.
\epr

\bprop\label{def-coinvolution-Mr}
On $M_r$, the formula $S_r(Qa)=Q(Sa)$, $a\in D(S)$, defines a proper coinvolution with the domain $D(S_r) = Q\big( D(S)\big)$.
\eprop
\bpr
Check first that $S_r$ is well-defined. If $Qa=0$ for $a\in D(S)$, then $a\in M_{**}^\perp$. Let us show first that $\widetilde S(a)\in M_{**}^\perp$. For $\mu\in M_{**}$,
$\mu(\widetilde S(a)) = \overline{ \mu^*(a)} =0$ since $\mu^*\in M_{**}$. But $S(a) = (S(a))^* \in M_{**}^\perp$ since the latter is a self-adjoint ideal.

The property \ref{def-coinvolution}(1) follows from the fact that $Q$ is surjective and $\sigma$-weak continuous; \ref{def-coinvolution}(2,3) from it being a $*$-homomorphism. For \ref{def-coinvolution}(4), take $x\in D(S)$ and $\mu,\nu\in [M_{**}]$. If $\mu\circ S_r$ and $\nu\circ S_r$ extend to normal functionals on $M_r$ then $\mu\circ S_r\circ Q=\mu\circ Q\circ S$ and $\nu\circ Q\circ S$ extend to normal functionals on $M$, and
\begin{align*}
\tilde\Delta (Qx)(\nu\circ S_r\otimes\mu\circ S_r)
& = (S_r Q\otimes S_r Q)\Delta(x)(\nu\otimes\mu)
 = (QS\otimes QS) \Delta(x)(\nu\otimes\mu)
 \\&= \Delta(x)(\nu\circ Q\circ S\otimes\mu\circ Q\circ S)
= \Delta S(x)(\mu\circ Q\otimes\nu\circ Q)
\\&=  \tilde\Delta Q S(x)(\mu\otimes\nu)
=  \tilde\Delta S_r(Qx)(\mu\otimes\nu).
\end{align*}
\epr

\bcor\label{cor-Mr}
$M_r$ with the quotient structure is a quantum semigroup with involution, with the coinvolution $S_r$; moreover, $(M_r)_{**}$ is dense in $(M_r)_*$, and $S_r$ is closeable.
\ecor

\bprop\label{Z*-hom}
The restriction of $Z_*:M_*\to (M_r)_*$ to $M_{**}$ is a *-homomorphism.
\eprop
\bpr
First note that $\mu(x)=\mu(\zeta x)$ for every $\mu\in M_{**}$, $x\in M$. 
For $\mu,\nu\in M_{**}$, $x\in M$ we have:
\begin{align*}
Z_*(\mu\nu)\big(Qx\big) &= (\mu\nu)\big(ZQx\big) = (\mu\nu)(\zeta \cdot x) = (\mu\nu)(x)
 = (\mu\otimes\nu)\big(\Delta(x)\big)
 = (\mu\otimes\nu)\big((\zeta\otimes\zeta)\Delta(x)\big)
\\& = (\mu\otimes\nu)\big((ZQ\otimes ZQ)(\Delta(x))\big)
 = (Z_*\mu\otimes Z_*\nu)\big(\bar\Delta(Qx)\big)
 = (Z_*\mu Z_*\nu)(Qx),
\end{align*}
which implies $Z_*(\mu\nu)=Z_*(\mu)Z_*(\nu)$.

Again, for $\mu\in M_{**}$, $x\in D(S)$
\begin{align*}
Z_*(\mu^*)(Qx) &= \mu^*(ZQx) = \mu^*(\zeta\cdot x) = \mu^*(x)
 = \overline{ \mu\big( (Sx)^*\big)}
 \\&= \overline{ \mu\big( \zeta\cdot (Sx)^*\big)}
= \overline{ \mu\big( ZQ( *Sx)\big)} =
\overline{ Z_*\mu\big(  (QSx)^*\big)} =
\overline{ Z_*\mu\big( \! *\tilde S(Qx)\big)} =
(Z_*\mu)^*  (Qx),
\end{align*}
whence $Z_*(\mu^*)=(Z_*(\mu))^*$.
\epr

On $M_*$, $Z_*$ might not be a homomorphism.

\section{Representations vs. their coefficients}

Theorem \ref{rep-by-coeffs} below is an analogue of \cite[1.4.2]{enock} and of \cite[Proposition 2.5]{kustermans}. Though the proofs in this section are very close to the cited publications, our situation is somewhat more general and we prefer to include the proofs.

\btm\label{rep-by-coeffs}
Suppose $M_{**}$ is dense in $M_*$. Let $I$ be a set, and let $x_{\a\beta}$: $\a,\beta\in I$ be elements in $D(S)\subset M$ which satisfy the relations $S(x_{\a\beta})=x_{\beta\a}^*$ and
\beq\label{Delta-for-pi-mult}
\Delta(x_{\a\beta}) = \sum_{\gamma} x_{\a\gamma}\otimes x_{\gamma\beta},
\eeq
the series converging absolutely $\sigma$-weakly, for all $\a,\beta\in I$. Then there exists a *-representation $\pi$ of $M_{**}$ on $\ell^2(I)$ such that
\beq\label{mu-pi}
\mu(x_{\a\beta}) = \langle \pi(\mu)e_\beta,e_\a\rangle \text{ for all } \mu\in M_{**} \text{ and all } \a,\beta\in I.
\eeq
\etm

Recall \cite[11.3]{palmer} that a positive functional $\phi$ on a Banach *-algebra $\cal A$ is called representable if there exists a *-representation $T: \cal A\to B(H)$ and a (topologically) cyclic vector $\xi\in H$ such that $\phi(a) = \langle T(a)\xi,\xi\rangle$ for all $a\in \cal A$. A positive functional is representable if and only if there is a constant $C$ such that $\phi(a)^2\le C\phi(a^*a)$, $a\in \cal A$. In this case $\|\phi\|=\|\xi\|^2\le C$. (See \cite[11.3.4]{palmer}.) Below we consider the algebra $M_{**}$ with its norm $\|\cdot\|_*$ which makes is a Banach $*$-algebra.

\blm
In the assumptions of Theorem \ref{rep-by-coeffs}, let $\xi\in\ell^2(I)$ be finitely supported, and let $\omega_\xi = \sum_{\a,\beta\in I} \xi_\a\bar\xi_\beta x_{\a\beta}$. Then $\omega_\xi$ is a positive representable functional on $M_{**}$ and $\|\omega_\xi\|_{M_{*\!*}^*}\le \|\xi\|_{\ell^2(I)}^2$.
\elm
\bpr
For $\mu\in M_{**}$,
$$
\omega_\xi(\mu^*\mu) = \Delta(\omega_\xi) (\mu^*\otimes\mu) = \sum_{\a,\beta} \xi_\a\bar\xi_\beta \Big( \sum_\gamma x_{\a\gamma}\otimes x_{\gamma\beta}\Big) (\mu^*\otimes\mu)
 = \sum_\gamma \Big( \sum_{\a,\beta} \xi_\a\bar\xi_\beta x_{\a\gamma} (\mu^*) x_{\gamma\beta}(\mu) \Big),
$$
where we change the summation order by absolute convergence. Since $x_{\a\gamma}(\mu^*) = (Sx_{\a\gamma})(\bar\mu) = (x_{\gamma\a}^*)(\bar\mu) = \overline{x_{\gamma\a}(\mu)}$,
$$
\omega_\xi(\mu^*\mu) = \sum_\gamma \Big( \sum_{\a,\beta} \xi_\a\bar\xi_\beta \overline{x_{\gamma\a}(\mu)} x_{\gamma\beta}(\mu) \Big)
= \sum_\gamma \Big| \sum_{\a} \bar\xi_\a  x_{\gamma\a}(\mu) \Big|^2 \ge 0,
$$
and we see that $\omega$ is positive.
Next,
$$
|\omega_\xi(\mu)|^2 = \Big| \sum_{\a} \xi_\a \sum_{\beta} \bar\xi_\beta x_{\a\beta}(\mu) \Big|^2
\le \|\xi\|_{\ell^2(I)}^2 \sum_{\a} \Big| \sum_{\beta} \bar\xi_\beta x_{\a\beta}(\mu) \Big|^2
 = \|\xi\|_{\ell^2(I)}^2 \,\omega_\xi(\mu^*\mu).
$$
This shows that $\omega_\xi$ is representable $\|\omega_\xi\|_{M_{*\!*}^*}\le \|\xi\|_{\ell^2(I)}^2$.
\epr

{\it Proof of the Theorem.\/}
Let $(e_\a)_{\a\in I}$ be the canonical base in $\ell^2(I)$. For $\mu\in M_{**}$ and finitely supported $\xi=\sum_a \xi_a e_\a\in \ell^2(I)$, set
$\pi(\mu)\xi = \sum_{\a,\beta} \xi_a x_{\beta\a}(\mu) e_\beta$. By the calculation in the Lemma above,
$$
\Big\|\sum_{\a,\beta} \xi_a x_{\beta\a}(\mu) e_\beta\Big\|^2
 = \sum_\beta \Big| \sum_\a \xi_a x_{\beta\a}(\mu) \Big|^2
 = \omega_{\,\bar{\!\xi}} (\mu^*\mu) \le \|\omega_{\,\bar{\!\xi}}\| \, \|\mu^*\|_* \|\mu\|_* = \|\xi\|^2 \|\mu\|_*^2.
$$
This shows that $\pi(\mu)$ extends to a bounded operator on $\ell^2(I)$ of norm at most $\|\mu\|_*$. By definition, $x_{\a\beta}(\mu) = \langle\pi(\mu)e_\beta,e_a\rangle$.

For $\mu,\nu\in M_{**}$,
\begin{align*}
\langle\pi(\mu\nu)e_\a,e_\beta\rangle &= x_{\beta\a}(\mu\nu) = \Delta(x_{\beta\a}) (\mu\nu) = \sum_\gamma x_{\beta\gamma}(\mu) x_{\gamma\a}(\nu)
= \sum_\gamma \langle \pi(\mu) e_\gamma, e_\beta\rangle  \langle \pi(\nu) e_a, e_\gamma\rangle
\\&= \langle \pi(\nu) e_a, \pi(\mu)^* e_\beta\rangle = \langle \pi(\mu)\pi(\nu) e_a, e_\beta\rangle.
\end{align*}
This proves that $\pi$ is multiplicative. Next,
$$
\langle \pi(\mu^*)e_\a,e_\beta\rangle = x_{\beta\a}(\mu^*) = \bar\mu(Sx_{\beta\a}) = \bar\mu(x_{\a\beta}^*) = \overline{\mu(x_{\a\beta})}
 = \overline{\langle \pi(\mu) e_\beta,e_\a\rangle} = \langle e_\a,\pi(\mu)e_\beta\rangle,
$$
so that $\pi$ is involutive on $M_{**}$.
\qed

\bcor\label{antirep-by-coeffs}
Suppose $M_{**}$ is dense in $M_*$. Let $I$ be a set, and let $x_{\a\beta}$: $\a,\beta\in I$ be elements in $D(S)\subset M$ which satisfy the relations $S(x_{\a\beta})=x_{\beta\a}^*$ and
\beq\label{Delta-for-pi-antimult}
\Delta(x_{\a\beta}) = \sum_{\gamma} x_{\gamma\beta}\otimes x_{\a\gamma},
\eeq
then there exists a *-antirepresentation $\pi$ of $M_{**}$ on $\ell^2(I)$ with the condition \eqref{mu-pi}.
\ecor
\bpr
Set $z_{\a\beta} = x_{\beta\a}$, then $(z_{\a\beta})_{\a,\beta\in I}$ satisfy the assumptions of Theorem \ref{rep-by-coeffs}. For the *-representation $\rho$ of $M_{**}$ on $\ell^2(I)$ given by the Theorem,
\[
\mu(x_{\a\beta}) = \mu(z_{\beta\a}) = \langle \rho(\mu)e_\a,e_\beta\rangle
\]
for all $\mu\in M_{**}$ and all $\a,\beta\in I$. Consider now the $\overline{\ell^2(I)}$ with the conjugate-linear structure and the product $\langle x,y\rangle_- = \langle y,x\rangle$. For every $\mu$, $\rho(\mu)^*$ is a linear operator on $\overline{\ell^2(I)}$, which we can denote $\pi(\mu)$. Moreover,
$$
\langle \pi(\mu)e_\beta,e_\a\rangle_- = \langle \rho(\mu)e_\a,e_\beta\rangle = \mu(x_{\a\beta}),
$$
so that $\pi$ satisfies \eqref{mu-pi}. Finally, $\pi$ is obviously anti-multiplicative and involutive.
\epr

\bprop\label{Delta-on-coeffs}
Let $\pi$ be a *-representation of $M_*$, and let $(x_{\a\beta})_{\a,\beta\in I}$ be its coefficients in a basis $(e_\a)_{\a\in I}$. Then the series \eqref{Delta-for-pi-mult} converges absolutely $\sigma$-weakly for all $\a,\beta\in I$.
\eprop
\bpr
This is a direct verification: for $\mu,\nu\in M_*$,
\begin{align*}
\Big(\sum_{\gamma} x_{\a\gamma}\otimes x_{\gamma\beta}\Big) (\mu\otimes\nu) &= \sum_{\gamma} x_{\a\gamma}(\mu) x_{\gamma\beta}(\nu)
 = \sum_{\gamma} \langle\pi(\mu)e_\gamma,e_\a\rangle \langle \pi(\nu)e_\beta,e_\gamma\rangle
\\& = \langle \pi(\nu)e_\beta, \pi(\mu)^*e_\a\rangle = \langle\pi(\mu)\pi(\nu)e_\beta,e_\a\rangle
= \langle\pi(\mu\nu)e_\beta,e_\a\rangle
= \Delta(x_{\a\beta})(\mu\otimes\nu).
\end{align*}
The series converges absolutely, since this is a decomposition of a scalar product.
\epr

\section{Absolutely continuous ideal}\label{sec-ideal}

Every *-representation of $M_{**}$ is bounded in the norm $\|\cdot\|_*$, but not necessarily in the norm of $M_*$. In \cite{kustermans}, J.~Kustermans proves that if $M$ is a locally compact quantum group then one implies the other. But in general we have no reason to expect this. Thus, we restrict the class of representations to those which have a ``generator'' (defined below) in $B(H)\bar\otimes M$, what implies in particular that they are completely bounded as representations of $M_*$ (but maybe not involutive, and maybe the image of $M_*$ is not self-adjoint).

There are some complications arising from the fact that $M_{**}$ might not be dense in $M_*$. Recall that this never occurs if $M$ is a locally compact quantum group \cite[Section 3]{kustermans}.

By Corollary \ref{cor-Mr}, we can always pass to $M_r$ and have $(M_r)_{**}$ dense in $(M_r)_*$.
In the sequel, we define the dual von Neumann algebra $\hat M$ first in the case when $M_{**}$ is dense, and set then $\hat M := \widehat{M_r}$. Next, we show that $\hat M$ has a structure of a quantum semigroup with involution.

\subsection{The case of $M_{**}$ dense in $M_*$}

Let $M$ be as above a quantum semigroup with involution, and suppose that $M_{**}\subset M_*$is dense. Let $\pi$ be an involutive  representation of $M_{**}$ on a Hilbert space. It is automatically continuous with respect to $\|\cdot\|_*$, but not necessarily with respect to $\|\cdot\|$. But in the case it is, $\pi$ extends by continuity to $M_*$, and this extension is a homomorphism (with a possibly non-selfadjoint range).
Representations of $M_*$ which are involutive on $M_{**}$ will be called below just {\it *-representations\/} of $M_*$.

$M_*$ can be considered as an operator space with the structure induced by duality to $M$.
Recall the following fact.

\bprop\label{cb=U}
Let $\pi$ be a representation of $M_*$ on a Hilbert space $K$. Then $\pi$ is completely bounded if and only if there exists $U\in B(K)\bar\otimes M\simeq (B(K)_*\hat\otimes_{\rm op}M_*)^*$ such that $U(\omega,\mu) = \omega(\pi(\mu))$ for every $\mu\in M_*$, $\omega\in N(K)$. In this case $\|U\|=\|\pi\|_{cb}$.
\eprop

\begin{notations}\label{nots-coeffs}
Let $\pi$ be a representation of $M_*$ on a Hilbert space $K$ with a basis $(f_\a)$. Denote by $\pi_{\a\beta}\in M$ the linear functional on $M_*$ defined by $\pi_{\a\beta}(\mu)=\langle \pi(\mu)f_\beta,f_\a\rangle$, $\mu\in M_*$. 
\end{notations}

In the next definition, the equivalence of conditions (1) and (2) has been proved in \cite[Theorem 3.5]{haar}. The proof makes no use of involution on $M_*$ or of the coinvolution (coinvolution) on $M$, so it remains valid in the present case. In fact, it is not even necessary to suppose that $\pi$ is a representation.

\begin{definition}
Let $\pi: M_*\to B(K)$ be a completely bounded $*$-representation of $M_*$ on a Hilbert space $K$.
Then $\pi$ is called {\it unitary} if one of the following equivalent conditions holds:
\begin{enumerate}
\item Exists a unitary $U\in B(K)\bar\otimes M$ called a {\it generator} of $\pi$ such that
\beq\label{U-pi}
U(\omega,\mu) = \omega(\pi(\mu))
\eeq
for every $\mu\in M_*$, $\omega\in B(K)_*$;
\item $\pi$ is non-degenerate and in some basis of $K$,
\beq\label{def-standard}
\sum_\gamma \pi_{\gamma\a}^*\cdot\pi_{\gamma\beta} = \sum_\gamma \pi_{\a\gamma}\cdot\pi_{\beta\gamma}^*
 = \begin{cases} 1, &\a=\beta\\ 0, &\a\ne\beta\end{cases}
\eeq
for every $\a,\beta$, the series converging absolutely in the $M_*$-weak topology of $M$.
\end{enumerate}
\end{definition}

From the Theorem 3.5 of \cite{haar} it follows that this definition does not depend, in fact, on the choice of a basis.

Moreover, one can verify that the following is true:
\bprop\label{construct-U-from-series}
Let the elements $\pi_{\a\beta}\in M$, $\a,\beta\in A$, be such that \eqref{def-standard} holds. Then there is a Hilbert space $K$ with a basis $(e_\a)_{\a\in A}$, a unitary operator $U\in B(K)\otimes M$ and a linear completely contractive map $\pi:M_*\to B(K)$ such that $\pi_{\a\beta}(\mu) = \langle\pi(\mu)e_\beta,e_\a\rangle = U(\omega_{\a\beta},\mu)$ for all $\mu\in M_*$ and $\a,\beta\in A$.
\eprop

For the proof, one realizes $M$ on a Hilbert space $H$ and defines $U$ as an operator on $K\otimes H$ by
$$
\langle U(e_\a\otimes x),(e_\beta\otimes y)\rangle = \pi_{\beta\a}( \mu_{x,y})
$$
for $x,y\in H$, $\a,\beta\in A$. It follows then from \eqref{def-standard} that $U$ is unitary, what implies the existence of $\pi$ with required properties.

Similarly to \cite{haar}, we give the following definitions:

\begin{definition}\label{def-Mtimes}
We call {\it non-unitary} the *-representations of $M_{**}$, completely bounded or not, which are not unitary. Let $M_{**}^\times\subset M_{**}$ be the (intersection of $M_{**}$ with the) common kernel of all {\it irreducible} non-unitary representations. If there are none, let $M_{**}^\times=M_{**}$. This is a two-sided *-ideal in $M_{**}$, which is called the {\it absolutely continuous ideal} of $M_*$. With the structure inherited from $M_{**}$, $M_{**}^\times$ is a Banach *-algebra.
\end{definition}

By definition, $M_{**}^\times$ is contained in $M_{**}$ and might not be an ideal in $M_*$. However, if $M_{**}$ is dense in $M_*$, one sees that the norm closure $[M_{**}^\times]$ is an ideal in $M_*$.

Note that the direct sum of a unitary and non-unitary representation is non-unitary; this is the reason to consider only irreducible representations in the definition above.

The main property of $M_{**}^\times$ is that every representation of $M_*$ which is irreducible on $M_{**}^\times$ must be unitary. For the duality construction we need more: that unitary representations do not vanish on it. To guarantee this, we exclude all degenerate cases by the following definition:

\begin{definition}\label{def-M0**}
Let $I^0$ be the weakly closed ideal in $M$ generated by $(M_{**}^\times)^\perp$, that is by the annihilator of $M_{**}^\times$. Set $M_{**}^0=M_{**}^\times$ if $I^0\ne M$ and $M_{**}^0=\{0\}$ otherwise.
\end{definition}

\brem\label{nondegen=contained-in-closure}
The following can be said also of any $*$-ideal of $M_{**}$. Every non-degenerate representation of $M_{**}^0$ extends uniquely to $M_{**}$ \cite[11.1.12]{palmer}. It is easy to show that $M_{**}$ is mapped into the weak closure of the image of $M_{**}^0$, i.e. to the von Neumann algebra generated by $M_{**}^0$. Conversely, if $\phi_1,\phi_2: M_{**}\to N$ are two *-homomorphisms to a von Neumann algebra $N$ which agree on $M_{**}^0$ and are such that $\phi_i(M_{**})$ is contained in the weak closure of $\phi_i(M_{**}^0)$ for $i=1,2$, then $\phi_1=\phi_2$.
\erem

\bprop\label{unitary-is-non-degenerate}
If $M_{**}^0\ne\{0\}$ then every unitary representation of $M_*$ is non-degenerate on $M_{**}^0$.
\eprop
\bpr
Let $\pi: M\to B(H)$ be a unitary representation of $M_*$. First, it is nonzero on $M_{**}^0$: otherwise we would have $\pi_{\a\beta}\in I^0$ for all its coefficients, and by \eqref{def-standard} this would imply $1\in I^0$, what is not true by assumption.

Next suppose that $\pi$ is degenerate on $M_{**}^0$. Let $L\subset H$ be the null subspace of $\pi(M_{**}^0)$. Then $\pi|^L$ is also unitary (with generator $(P_L\otimes 1)U$ where $P_L$ is the projection onto $L$ and $U$ is the generator of $\pi$) and vanishes on $M_{**}^0$, what is impossible.
\epr

\bprop\label{nondegenerate-is-cb}
Every $*$-representation of $M_{**}^0$ extends to a completely contractive representation of $M_*$.
\eprop
\bpr
Every irreducible representation of $M^0_{**}$ is by definition unitary and thus extends to a completely contractive representation of $M_*$. In general, the norm of a *-representation $\pi: M^0_{**}\to B(H)$ is majorated by the supremum over all irreducible representations; thus $\|\pi(x)\|\le \|x\|_{M_*}$ for every $x\in M^0_{**}$, and it follows that $\pi$ extends to a contractive representation of $M_*$ by continuity.

To show that $\pi$ is completely bounded, we need the disintegration theory in its a priori non-separable form \cite{henrichs}. Let $A$ be the unital $C^*$-algebra generated by $\pi(M^0_{**})$ and $\Id_H$. There exists a measure $\nu$ on the state space $\Omega$ of $A$ such that the identity representation of $A$ on $H$ is isomorphic to the direct integral $\int_\Omega \rho_\phi d\nu(\phi)$ of irreducible representations $\rho_\phi$ of $A$ on respective spaces $H_\phi$; note that the direct integral is understood in the sense of W. Wils \cite{wils}, a definition suitable for non-separable disintegration.

Every $\rho_\phi\circ\pi$ is either irreducible on $M^0_{**}$ or vanishes on it. This allows to verify explicitly that the matrix norm of $\pi_n: M_n(M^0_{**})\to M_n(B(\int_\Omega H_\phi))$ is majorated by 1, what proves that $\pi$ extends to a completely contractive representation of $M_*$.
\epr

The following is a generalization of Theorem 5.5 of \cite{haar}. The proof is changed in significant details to treat possible inequality of $M_*$ and $M_{**}$.

\btm\label{nondegenerate-is-standard}
Every bounded $*$-representation $\pi$ of $M_*$ which is nondegenerate on $M_{**}^0$ is unitary.
\etm
\bpr
By Proposition \ref{nondegenerate-is-cb} (extension is the same as $\pi$ by non-degeneracy), $\pi$ is in fact completely contractive.
By Proposition \ref{cb=U}, there exists $U\in B(H)\bar\otimes M$ such that $U(\omega,\mu) = \omega(\pi(\mu))$, and all we need is to prove that $U$ is unitary.

Let $M$ be realized on a Hilbert space $K$. We will need several times the following representation. Fix $x\in K$, $\xi\in H$ and bases $(e_\a)\subset K$, $(f_\beta)\subset H$. Then
$$
\langle U(\xi\otimes x),f_\beta\otimes e_\a\rangle = \langle \pi(\mu_{x e_\a}) \,\xi,f_\beta\rangle,
$$
so that
\begin{align}\label{U-series}
U(\xi\otimes x)
 &= \sum_{\a,\beta} \langle U(\xi\otimes x),f_\beta\otimes e_\a\rangle\, f_\beta\otimes e_\a
 = \sum_{\a,\beta} \langle \pi(\mu_{x e_\a}) \,\xi,f_\beta\rangle \,f_\beta\otimes e_\a                      \notag
\\ &= \sum_\a \Big( \sum_\beta \langle \pi(\mu_{x e_\a})\, \xi,f_\beta\rangle \;f_\beta \Big) \otimes e_\a
= \sum_\a \pi(\mu_{x e_\a})\, \xi \otimes e_\a
\end{align}
(convergence is everywhere in the Hilbert space norm). Similarly,

\begin{align}\label{Ustar-series}
U^*(\xi\otimes x)
 &= \sum_{\a,\beta} \langle \xi\otimes x,U(f_\beta\otimes e_\a)\rangle\, f_\beta\otimes e_\a
 = \sum_{\a,\beta} \langle \pi(\mu_{e_\a x}) \,f_\beta, \xi\rangle^- \,f_\beta\otimes e_\a                      \notag
\\ &= \sum_\a \Big( \sum_\beta \langle \pi(\mu_{e_\a x})^*\, \xi,f_\beta\rangle \;f_\beta \Big) \otimes e_\a
= \sum_\a \pi(\mu_{e_\a x})^*\, \xi\otimes e_\a
\end{align}

We can suppose that $M$ is realized in its standard form. Then every $\mu\in M_*$ is equal to $\mu_{xy}$ for some $x,y\in K$, so that we can identify $M_*$ with $K\otimes \bar K$.

For subspaces $E,F\subset K$, let $M_{E,F}\subset M_*$ denote the closed subalgebra generated by $\mu_{xy}$ with $x\in E$, $y\in F$. Denote also $M^0_{E,F} = M_{E,F}\cap M^0_{**}$. These subalgebras are not supposed to be self-adjoint. By $M_{E,F}^\times$, $M_{E,F}^{0\times}$ we denote the $\|\cdot\|_*$-closed $*$-subalgebras generated by $M_{E,F}\cap M_{**}$ and $M_{E,F}^0$ respectively. Note that by definition $M^0_{E,F}\subset M_{E,F}^{0\times}$ but it might be $M_{E,F}\not\subset M^\times_{E,F}$.

\blm\label{invar-subspace-U}
A closed subspace $L\otimes E\subset H\otimes K$ is $U$-invariant if and only if $\pi(M_{E,K})L\subset L$ and $\pi(M_{E,E^\perp})L=\{0\}$. It is $U^*$-invariant if and only if $\pi(M_{K,E})^*L\subset L$ and $\pi(M_{E^\perp,E})^*L=\{0\}$.
\elm
The proof is identical to the proof of Lemma 5.6 of \cite{haar}.

\blm\label{invar-subspace-F-UF}
For every separable subspace $V=G\otimes F\subset H\otimes K$ there exist closed separable subspaces $E\subset K$, $L\subset H$ such that $V\cup UV \cup U^*V\subset L\otimes E$ and $M_{F,F}\subset \overline{M_{**}\cap M_{E,E}}$.
\elm
\bpr
It follows from the assumptions that $F$ and $G$ are separable. Let $X\subset M_{**}$ be countable and such that $M_{F,F}\subset \overline{X}$. By the isomorphism $M_*\simeq K\otimes\bar K$, we can write $X=\{\mu_{x_ny_n}:n\in\N\}$. Set $F_1={\rm lin}\{x_n,y_n: n\in\N\}$, then $X\subset M_{F_1,F_1}$.

Since $V$ is separable, so is $V_2=\overline{G\otimes F_1+V+UV+U^*V}$. Pick a sequence $(v_n)$ dense in $V$ and orthonormal bases $(e_\a)\subset K$, $(f_\beta)\subset H$. Every $v_n$ is contained in $\overline{\rm lin}\{ f_\beta\otimes e_\a: \a\in A_n, \beta\in B_n\}$ with countable $A_n$, $B_n$. Then for $E=\overline{\rm lin}\{ e_\a: \a\in \cup A_n\}$ and $L = \overline{\rm lin}\{ f_\beta: \beta\in \cup B_n\}$ we have $V\cup UV \cup U^*V\subset L\otimes E$, and $X\subset M_{E,E}$ so that $E,L$ are as required.
\epr

\blm
Every $v\in K\otimes H$ can be embedded into a $U,U^*$-invariant separable subspace $L\otimes E$ such that 
$L$ is essential for $\pi(M_{E,K}^0)|_L$ and $M_{E,E}\subset \overline{M_{**}\cap M_{E,E}}$.
\elm
\bpr
Construct separable subspaces $E_k$, $L_k$ by induction as follows. Let $L_1\otimes E_1$ be any separable subspace containing $v$. Suppose now that $E_{k-1}$, $L_{k-1}$ are constructed for some $k\ge2$. Since $H$ is essential for $M_{**}^0$, there are sequences $\mu^{(k)}_n\in M_{**}^0$, $\xi^{(k)}_n\in H$ such that ${\rm lin}\{ \pi(\mu^{(k)}_n)\xi^{(k)}_n\}$ is dense in $L_{k-1}$. Since $M$ is in the standard form, every $\mu^{(k)}_n$ can be represented as $\mu^{(k)}_n=\mu_{x^{(k)}_n,y^{(k)}_n}$ with $x^{(k)}_n,y^{(k)}_n\in K$. Set $E'_k=\overline{ E_{k-1}+{\rm lin}\{ x^{(k)}_n: n\in\N \}}$, $L'_k=\overline{L_{k-1}+ {\rm lin}\{ \xi^{(k)}_n: n\in\N \}}$. Then $\mu^{(k)}_n \in M_{E'_k,K}^0$ and $\xi^{(k)}_n\in L'_k$ for all $n$. By Lemma \ref{invar-subspace-F-UF} there are separable subspaces $E_k$, $L_k$ such that
$$L'_k\otimes E'_k\cup U(L'_k\otimes E'_k) \cup U^*(L'_k\otimes E'_k) \subset L_k\otimes E_k
$$
and $M_{E'_k,E'_k}\subset \overline{M_{**}\cap M_{E_k,E_k}}$.

Set  $E=\overline{\cup E_k}$ and $L=\overline{\cup L_k}$, then $L\otimes E = \overline{ \cup (L_k\otimes E_k)}$ since $E_k$, $L_k$ are increasing. We have $U(L_k\otimes E_k)\subset L_{k+1}\otimes E_{k+1}$ and $U^*(L_k\otimes E_k)\subset L_{k+1}\otimes E_{k+1}$
 for all $k$, what implies the $U,U^*$-invariance of $L\otimes E$. Moreover, by construction the set $\{\pi(\mu)\xi: \mu\in M^0_{E,K},\xi\in L\}$ is dense in $L$. Finally, $M_{E_k,E_k} \subset \overline{M_{**}\cap M_{E_{k+1},E_{k+1}}} \subset \overline{M_{**}\cap M_{E,E}}$ for every $k$, so $M_{E,E} = \overline{\cup M_{E_k,E_k}} \subset \overline{M_{**}\cap M_{E,E}}$.
\epr

{\it Proof of the theorem.} Take any $v\in H\otimes K$. Let $L\otimes E\subset H\otimes K$ be $U, U^*$-invariant and separable, such that $v\in L\otimes E$, $M_{E,E}\subset \overline{M_{**}\cap M_{E,E}}$ and $L$ is essential for $M_{E,K}^0$.
It follows that $\pi(M_{E,K})L\subset L$; $\pi(M_{E,E^\perp})L=\{0\}$; $\pi(M_{K,E})^*L\subset L$ and as a consequence $\pi(M_{E,E}^\times)L\subset L$.

Fix an orthonormal base $(e_\a)_{\a\in A}$ in $K$ such that the (countable) subset $(e_\a)_{\a\in A_1}$ is a base for $E$. Note that $(L\otimes E)^\perp$ is also $U,U^*$-invariant.
For $x\in E^\perp$, $\xi\in L$ one has $U^*(\xi\otimes x)\in (L\otimes E)^\perp$; from \eqref{Ustar-series} it follows that $\pi(\mu_{e_\a x})^*\xi\in L^\perp$ if $e_\a\in E$. It follows that $\pi(M_{E,E^\perp})^*L\subset L^\perp$.
For $x\in K$, $\xi\in L^\perp$ we have again $U^*(\xi\otimes x)\in (L\otimes E)^\perp$; and from \eqref{Ustar-series} it follows that $\pi(\mu_{e_\a x})^*\xi\in L^\perp$ if $e_\a\in E$. It follows that $\pi(M_{E,K})^*L^\perp\subset L^\perp$.

Let $C^*(X)$ denote the closed $*$-algebra generated by $X$ in its relevant space of operators, and let $r_L:B(H)\to B(L)$ be the reduction onto $L$. The reasoning above shows that $r_L: C^*(\pi(M_{E,K})) \to B(L)$ is a $*$-representation which vanishes on $\pi(M_{E,E^\perp})$. Denote $\rho=r_L\circ\pi$; we have $\rho(M_{E,E^\perp})=0$.

Recall that $M_{E,E}^\times$ is closed in the $\|\cdot\|_*$-norm; the assumption $M_{E,E}\subset \overline{M_{**}\cap M_{E,E}}$ implies moreover that in the norm of $M_*$ we have $M_{E,E} \subset \overline{M_{E,E}^\times}$, so that $\rho(M_{E,K}) = \rho(M_{E,E}) \subset \overline{\rho(M_{E,E}^\times)}=C^*(\rho(M_{E,E}^\times))$. From the other side, obviously $C^*(\rho(M_{E,E}^\times))\subset C^*(\rho(M_{E,E}))$, so in fact $C^*(\rho(M_{E,E}^\times))=C^*(\rho(M_{E,K}))$.

Set $\cal A=C^*(\rho(M_{E,K}^{0}))$. Since $M_{E,E}$ is separable, so are $\rho(M_{E,E})=\rho(M_{E,K})\supset \rho(M_{E,K}^{0})$ and $\cal A$. The identity representation of $\cal A$ in $B(L)$ is decomposed into a direct integral of irreducible representations \cite[8.5.2]{dixmier}: there exist a set $P$ equipped with a probability measure $\beta$; an integrable field of Hilbert spaces $\Gamma\subset \{(H_p)_{p\in P}\}$; a field of representations $\sigma_p: {\cal A}\to B(H_p)$, $p\in P$, where every $\sigma_p$ is irreducible; and an isometric isomorphism $V:L\to \Gamma=\int_P^\oplus H_p d\beta(p)$ such that: if $\xi\in L$ and $V(\xi)=\int_P \xi_p \,d\beta(p)$,
then for every $a\in \cal A$ we have $V(a\xi)=\int_P \sigma_p(a)\xi_p \,d\beta(p)$.

For every $\mu\in M_{E,K}^{0\times}$ and $\xi\in L$, denoting $V\xi=\int\xi_p$, we have
\beq\label{rho-equal-int}
V(\rho(\mu)\xi)=\int_P \rho_p(\mu)\xi_p \,d\beta(p)
\eeq
with irreducible representations $\rho_p=\sigma_p\circ\rho$ of $M^{0\times}_{E,K}$.

$M^{0\times}_{E,K}$ is an ideal in $M^\times_{E,K}$ as well as $\rho(M^{0\times}_{E,K})$ in $\rho(M^\times_{E,K})$, so we can extend $\sigma_p$ uniquely and irreducibly to $\rho(M_{E,K}^\times)$, still denoting this extension by $\sigma_p$. By continuity, $\sigma_p$ extends to $C^*(\rho(M_{E,K}^\times))=C^*(\rho(M_{E,K}))$, and this allows to extend $\rho_p$ to $M_{E,K}$.
These considerations show also that 
\eqref{rho-equal-int} is valid in fact for all $\mu\in M_{E,K}$.

For every $p$, $\sigma_p\circ r_L$ is lifted \cite[2.10]{dixmier}
from $C^*(\rho(M_{E,K}))$ to an irreducible representation $\tilde\sigma_p$ of $C^*(\pi(M_*))$, probably on a bigger space $\tilde H_p\supset H_p$. Set $\tilde\rho_p = \tilde \sigma_p\circ \pi$.

Since $\pi$ is completely bounded, every $\tilde\rho_p$ is completely bounded too. For $\mu\in M_{E,K}$, $\tilde\rho_p(\mu)|_{H_p}=\sigma_p\circ\rho(\mu)=\rho_p(\mu)$ .
In particular, $H_p$ is invariant under $\tilde\rho_p(M_{E,K})$,
and $\tilde\rho_p(M_{E,E^\perp})|_{H_p}=0$ (since $\rho(M_{E,E^\perp})=0$).
Being irreducible and nonzero on $M_{**}^0$, $\tilde\rho_p$ is unitary, with a unitary generator $U_p\in B(\tilde H_p)\otimes M$. By the reasoning above, we can apply Lemma \ref{invar-subspace-U} and conclude that $H_p\otimes E$ is $U_p$-invariant.

The rest of the proof is identical to that of \cite[Theorem 5.5]{haar}.
$\int(H_p\otimes E)_{p\in P}$ is also a field of Hilbert spaces, isomorphic to $L\otimes E$ under the isomorphism $\tilde V=V\otimes \id$. We have then for $x\in E$, $\xi\in L$ that $\mu_{x e_\a}\in M_{E,K}$ and so $\pi(\mu_{x e_\a})\xi=\rho(\mu_{x e_\a})\xi$. Denote $V\xi=\int \xi_p$; we have
\begin{align}
\tilde V(U(\xi\otimes x)) &= \tilde V\Big(\sum_{\a\in A_1} \pi(\mu_{x e_\a})\, \xi\otimes e_\a\Big)         \notag
\\&= \sum_{\a\in A_1} V(\rho(\mu_{x e_\a})\, \xi)\otimes e_\a
 = \sum_{\a\in A_1} \int \big(\rho_p(\mu_{x e_\a})\, \xi_p\big) \otimes e_\a                                \notag
\\&= \sum_{\a\in A_1} \int \big(\rho_p(\mu_{x e_\a})\, \xi_p\big)\otimes e_\a.                              \label{series-ea-int}
\end{align}
The last series converges in the Hilbert norm of $\int(H_p\otimes E)$.

From the other hand, for every $p$ we have a formula similar to \eqref{U-series}: if $x\in E$, $\xi_p\in H_p$, then
$$
U_p(\xi_p\otimes x) = \sum_{\a\in A_1} \rho_p(\mu_{x e_\a})\, \xi_p\otimes e_\a,
$$
so the series in \eqref{series-ea-int} converges pointwise to $\int U_p(\xi_p\otimes x)$. Both imply convergence in measure in the following sense: denote  $\phi_{\a p}=\rho_p(\mu_{x e_\a})\, \xi_p\otimes e_\a$, then for every $\e>0$
\begin{align*}
&\beta\{ p: \|\big(\tilde V(U(\xi\otimes x))\big)_p - \sum_{\a\in B} \int \phi_{\a p}\|\ge\e\} \to0,\\
&\beta\{ p: \|U_p(\xi_p\otimes x) - \sum_{\a\in B} \int \phi_{\a p}\|\ge\e\} \to0
\end{align*}
as finite set of indices $B\subset A_1$ increases (the reasoning for real-valued functions applies verbatim). It follows that $\big(\tilde V(U(\xi\otimes x))\big)_p = U_p(\xi_p\otimes x)$ almost everywhere, that is
$$
\tilde V(U(\xi\otimes x)) = \int U_p(\xi_p\otimes x) = \Big(\int U_p\Big)(\tilde V(\xi\otimes x)).
$$
It follows that $(U_p)$, or strictly speaking $(U_p|_{H_p\otimes E})$, is a measurable field of operators on $\int H_p\otimes E$, and $\tilde V U = \int U_p\tilde V$. Since $\int U_p$ is unitary, so is $U$ on $L\otimes E$.

As the initial  vector $v$ was arbitrary, we get that $U$ is unitary on $H\otimes K$, what proves the theorem.
\epr

\begin{definition}\label{def-Mhat}
Denote by $C^*(M_{**}^0)$ the \cst-enveloping algebra of $M_{**}^0$, that is the completion of $M_{**}^0$ with respect to the maximal \cst-seminorm on it. Denote by $\Mhat$ the enveloping von Neumann algebra of $C^*(M_{**}^0)$ and by $\Phi: M_{**}^0\to\Mhat$ the canonical map into it. As every non-degenerate representation, $\Phi$ extends uniquely to $M_{**}$ (see remark \ref{nondegen=contained-in-closure}), and by continuity to $M_*$. In the sequel, we always consider $\Phi$ as a map from $M_*$ to $\Mhat$. Note that a priori $\Phi$ might not be injective.\\
By Proposition \ref{nondegenerate-is-cb}, $\Phi$ is completely contractive, and by Theorem \ref{nondegenerate-is-standard} it is unitary, unless $\Mhat=\{0\}$. Denote also by $\hat\Phi$ the preadjoint map $\Phi_*: \Mhat_*\to M$.
\end{definition}

\brem
This definition coincides with the one given in \cite{haar} in the 
 case of a bounded coinvolution $S$. In this case, in particular, $M_{**}$ is always dense in $M_*$.
\erem

\subsection{The case of $M_{**}$ not dense in $M_*$}

\begin{definition}\label{def-Mhat-not-dense}
In general, set $\Mhat = \widehat{M_r}$, $\Phi_M = \Phi_{M_r}\circ Z_*$, where $Z_*:M_*\to M_{r*}$ is the map defined in Remark \ref{def-Z}.
By Proposition \ref{Z*-hom}, $\Phi_M$ is a *-homomorphism on $M_{**}$. However, if $M_{**}$ is not dense in $M_*$, then $\Phi_M$ might not be multiplicative on $M_*$. Still, $\Phi_M$ has a generator $V_M=(\Id\otimes Z)(V_{M_r})\in \Mhat\otimes M$, but it is not unitary, since $V_M$ is contained in the ideal $\Mhat\otimes (\zeta M)$ (see notations in Definition \ref{def-Z}).
\end{definition}

\bprop\label{M*-hom-to-vna}
Let $N$ be a von Neumann algebra and $\pi: M_*\to N$ a *-homomorphism. Then there is a unique *-homomorphism $\tilde\pi: \Mhat\to N$ such that $\tilde\pi\circ\Phi=\pi$ on $[M_{**}]$ and $\pi(M_*)$ is contained in the weak closure of $\pi(M^0_{**})$ in $N$. Is $M_{**}$ is dense in $M_*$, then the equality holds on $M_*$.
\eprop
\bpr
From $M^0_{**}$, $\pi$ is lifted by the universality property to $\tilde\pi: \Mhat\to N$, so that $\tilde\pi\circ\Phi=\pi$ (on $M^0_{**}$). Now $\pi$ and $\tilde\pi\circ\Phi$ are two extensions of $\pi$ to $M_{**}$ with the property that the image of $M_{**}$ is contained in the weak closure of $\pi(M^0_{**})$. As for any *-ideal in $M_{**}$, this implies that $\pi=\tilde\pi\circ\Phi$ on $M_{**}$. Finally, by density this equality holds also on $[M_{**}]$, and on the whole of $M_*$ is $M_{**}$ is dense in it.
\epr

\section{The dual Hopf-von Neumann algebra}\label{section-dual-von-neumann}

In this section we still suppose that $M_{**}$ is dense in $M_*$. General case will follow for granted since it is reduced to the case of $M_{**}$ dense.

In order to define a comultiplication on $\hat M$, we define first the Kronecker product $\Phi\times\Phi$ of $\Phi$ by itself, via its coordinates. Recall the procedure of how it is done. If $\hat M\subset B(K)$ and $(e_\a)_{\a\in A}$ is a base of $K$, let $\Phi_{\a\beta}$ be the coefficients of $\Phi$. Then $\Phi_{\a,\beta}\cdot \Phi_{\a',\beta'}$ with $\a,\beta,\a',\beta'\in A$ satisfy the equalities \eqref{def-standard}, so they generate a unitary $V\in \hat M\otimes\hat M\otimes M$ such that $V(\omega_{\a\beta},\omega_{\a',\beta'},\mu) = (\Phi_{\a,\beta}\cdot \Phi_{\a',\beta'})(\mu)$ for all $\mu\in M_*$. The corresponding map $\Phi\times\Phi:M_*\to B(K\otimes K)$ has clearly its range in $\hat M\otimes\hat M$ and is multiplicative on $M_*$, since $\Delta(\Phi_{\a,\beta}\cdot \Phi_{\a',\beta'}) = \sum_{\gamma,\gamma'} \Phi_{\a,\gamma}\Phi_{\a'\gamma'}\otimes\Phi_{\gamma\beta}\Phi_{\gamma',\beta'}$ (and this follows from the multiplicativity of $\Phi$).

\bprop\label{involutive=coeff_in_D(S)}
Let $\pi$ be a representation of $M_*$ in a Hilbert space $H$. Assume that $S$ is extended according to Corollary \ref{S-extended}. Then $\pi$ is involutive if and only if\/ for every $\omega\in B(H)_*$, $\pi_*(\omega)\in D(S)$ and $S(\pi_*(\omega))=\pi_*(\bar\omega)^*$.
\eprop
\bpr
Let $\pi$ be involutive. For $\mu\in M_{**}$,
$$
\pi_{\a\beta}(\mu^*) = \langle \pi(\mu^*)e_\beta,e_\a\rangle = \overline{ \langle \pi(\mu)e_\a,e_\beta\rangle}
 = \overline{\pi_{\beta\a}(\mu)}.
$$
By definition \ref{def-R}, this implies $\pi_{\a\beta}\in D(R)$ and $R(\pi_{\a\beta})=\pi_{\beta\a}$, so that $S(\pi_*(\omega_{\a\beta}))=\pi_*(\bar\omega_{\a\beta})^*$. By the weak continuity of $R$, the statement follows for every $\omega$.

Conversely, $R(\pi_{\a\beta})=\pi_{\beta\a}$ for all $\a,\beta$ implies by the same calculation that $\pi(\mu^*)=\pi(\mu)^*$ for every ${\mu\in M_{**}}$.
\epr

\bprop
$\Phi\times\Phi$ is involutive on $M_{**}$.
\eprop
\bpr
Suppose first that $M_{**}$ is dense in $M_*$. For $\mu\in M_{**}$ and $\omega,\upsilon\in\Mhat_*$,
$(\Phi\times\Phi)_*(\omega\otimes\upsilon) = \hat\Phi(\omega)\hat\Phi(\upsilon)\in D(S)$ since $\hat\Phi(\Mhat_*)\subset D(S)$, and
$$
S\big((\Phi\times\Phi)_*(\omega\otimes\upsilon)\big) = S\hat\Phi(\upsilon) S\hat\Phi(\omega)
= \hat\Phi(\bar\upsilon)^*\hat\Phi(\bar\omega)^* = \big( (\Phi\times\Phi)_*(\bar\omega\otimes\bar\upsilon)\big)^*.
$$
By Proposition \ref{involutive=coeff_in_D(S)}, $\Phi\times\Phi$ is involutive.

We have $(\Phi\times\Phi)_{\a,\a',\beta,\beta'} = \Phi_{\a,\beta}\cdot \Phi_{\a',\beta'}\in D(R)$ for all $\a,\a',\beta,\beta'$, moreover $R((\Phi\times\Phi)_{\a,\a',\beta,\beta'}) = R(\Phi_{\a,\beta}) R(\Phi_{\a',\beta'}) = \Phi_{\a',\beta'}\cdot \Phi_{\a,\beta} = (\Phi\times\Phi)_{\a',\a,\beta',\beta}$. This proves the statement.
\epr

\bprop
$\Mhat_{**} \supset \hat\Phi^{-1}\big( \hat\Phi(\Mhat_*)\cap \hat\Phi(\Mhat_*)^*\big)$.
\eprop

By universality, $\Phi\times\Phi$ lifts to a *-homomorphism $\hat\Delta:\hat M\to \Mhat\otimes\Mhat$, so that $\hat\Delta\Phi = \Phi\times\Phi$ on $M^0_{**}$. Since both are non-degenerate, their extension to $M_{**}$ is unique, thus the equality holds on $M_{**}$; by continuity, it holds on $M_*$ as well.

\bprop\label{dual-comultiplication}
$\hat\Delta$ is a comultiplication on $\hat M$.
\eprop
\bpr
$\hat\Delta$ is involutive by construction. The fact that it is coassociative and unital is proved exactly as in Proposition 6.6 \cite{haar}.
\epr

Comultiplication on $\Mhat$ turns $\Mhat_*$, as usual, into a Banach *-algebra.

\bprop\label{hat-Phi-mult}
$\hat\Phi$ is multiplicative.
\eprop
\bpr
For $\upsilon,\omega\in\Mhat_*$, $\mu\in M_*$ by definition
$$
(\upsilon\cdot\omega)\big(\Phi(\mu)\big) = (\upsilon\otimes\omega)\big(\hat\Delta\Phi(\mu)\big)
 = (\upsilon\otimes\omega)\big((\Phi\times\Phi)(\mu)\big).
$$
At the same time, $(\upsilon\cdot\omega)\big(\Phi(\mu)\big)=\hat\Phi(\upsilon\cdot\omega)(\mu)$ and $(\upsilon\otimes\omega)\big((\Phi\times\Phi)(\mu)\big) = \big(\hat\Phi(\upsilon)\hat\Phi(\omega)\big)(\mu)$. This implies that $\hat\Phi(\upsilon\cdot\omega) = \hat\Phi(\upsilon)\hat\Phi(\omega)$.
\epr

\begin{definition}\label{def-hat-S}
Let $V\in \Mhat\otimes M$ be the unitary generator of $\Phi$. On $D(\hat S) = \Phi(M_*)$, define a map $\hat S$  
by $(\hat S\Phi\mu, \omega) = V^*(\omega,\mu)$, $\mu\in M_*$, $\omega\in \hat M_*$.
\end{definition}

\bprop\label{S-Phi-mu}
If $\mu\in M_*$ is such that $\mu\circ S$ extends to a normal functional on $M$, then
$\hat S\Phi(\mu)(\omega) = (\mu\circ S)\big(\hat\Phi(\omega)\big)$ for all $\omega\in \Mhat_*$.
\eprop
\bpr
For any $\mu\in M_*$ and $\omega\in \Mhat_*$, we have:
$$
\hat S\Phi(\mu)(\omega) = V^*(\omega,\mu) = \overline{V(\bar\omega,\bar\mu)}
 = \bar\omega\big( \Phi(\bar\mu)\big) = \omega\big( \Phi(\bar\mu)^*\big),
$$
so that
\beq\label{hatSPhi}
\hat S\Phi(\mu) = \Phi(\bar\mu)^*.
\eeq
The condition on $\mu$ is equivalent to the fact that $\bar\mu\in M_{**}$; if it holds, we can continue the calculation above as
$$
\hat S\Phi(\mu)(\omega) = \bar\mu^*\big( \hat\Phi(\omega) \big) = (\mu\circ S)\big(\hat\Phi(\omega)\big),
$$
which proves the statement.
\epr

\bprop\label{hatPhi(Mhat)-in-D(S)}
If $M_{**}$ is dense in $M_*$, then $\hat\Phi(\Mhat_*)\subset D(S)$ and $S\hat\Phi(\omega) = \hat\Phi(\bar\omega)^*$.
\eprop
\bpr
Follows immediately from Proposition \ref{involutive=coeff_in_D(S)}.
\epr

Suppose that $M\subset B(K)$ and denote by $\theta$ the flip map on $K\otimes K$. Let $H$ be another Hilbert space. Recall the leg numbering notation: for $V\in B(H)\otimes M$, we denote $V_{12}=V\otimes\Id$ and $V_{13}=(\Id\otimes\theta)(V\otimes \Id)(\Id\otimes\theta)$. Exactly as in \cite[1.5.1]{enock}, one can prove:

\begin{prop}\label{prop-mult-antimult}
Let $\pi : M_*\to B(H)$ be a completely bounded linear map, and let $V\in B(H)\otimes M$ be such that $V(\omega,\mu)=\omega(\pi(\mu))$ for all $\omega\in B(H)_*$, $\mu\in M_*$. Then
\begin{enumerate}
\item $\pi$ is multiplicative if and only if
\beq\label{V-mult}
(\Id\otimes\Delta)(V)
= V_{12}V_{13};
\eeq
\item $\pi$ is anti-multiplicative if and only if
\beq\label{V-antimult}
(\Id\otimes\Delta)(V) 
= V_{13} V_{12}.
\eeq
\end{enumerate}
\end{prop}

\begin{prop}\label{dual-coinvolution}
$\hat S$ is a proper coinvolution on $\hat M$.
\end{prop}
\bpr
(1) $D(\hat S)=\Phi(M_*)$ is by definition dense in $\Mhat$.

(2) $D(\hat S)$ is closed under multiplication since $\Phi$ is a homomorphism.
 Let $V$ be the generator of $\Phi$. Since $\Phi$ is a homomorphism, by Proposition \ref{prop-mult-antimult}
$$
(\Id\otimes\Delta)(V)
= V_{12}V_{13}.
$$
For $V^*$ we have, since $\Delta$ is involutive and $\theta^*=\theta$:
$$
(\Id\otimes\Delta)(V^*)
 = (V^*)_{13} (V^*)_{12}.
$$
By definition, $V^*$ is the generator of $\hat S\Phi$, thus by the same Proposition $\hat S\Phi$ is anti-multiplicative. It follows that $\hat S$ is anti-multiplicative on its domain.

(3) A direct calculation shows that
$$
V^*(\bar\omega,\mu) = \overline{V(\omega,\bar\mu)}
$$
for all $\mu\in M_*$, $\omega\in \Mhat_*$.

Next,
\beq\label{omega-involution}
\omega\big(\Phi(\bar\mu)^*\big) = \overline{ \bar\omega\big( \Phi(\bar\mu)\big) }
 = \overline{ V(\bar\omega,\bar\mu)} = V^*(\omega,\mu)
  = \omega\big( \hat S\Phi(\mu)\big)
\eeq
It follows that $*\hat S\Phi\mu = \Phi(\bar\mu)$ for every $\mu\in M_*$. Now we get immediately $*\hat S(*\hat S\Phi\mu) = *\hat S(\Phi\bar\mu) = \Phi(\bar{\bar\mu}) =\Phi(\mu)$ and as a consequence $(*\hat S)^2=\Id$ as required.

(4) 
Let $u,v\in \Mhat_*$ be such that $u\circ \hat S$ and $v\circ \hat S$ extend to normal functionals on $\Mhat$.
We need to show that for every $\mu\in M_*$,
\beq\label{DeltaHat}
\hat \Delta \big(\Phi(\mu)\big)(v\circ \hat S\otimes u\circ \hat S) = \hat \Delta \hat S\Phi(\mu)(u\otimes v).
\eeq
The left hand side equals to
\beq\label{DeltaPhi1}
\hat \Delta (\Phi\mu)(v\circ \hat S\otimes u\circ \hat S)
= (\Phi\times\Phi)(\mu) (v\circ \hat S\otimes u\circ \hat S)
= \mu\big( \hat\Phi(v\circ \hat S)\cdot\hat\Phi(u\circ \hat S)\big).
\eeq
One checks that for any $\nu\in M_*$
$$
\nu(\hat\Phi(u\circ \hat S)) = \hat S\Phi(\nu)(u) = V^*(u,\nu) = \overline{V(\bar u,\bar \nu)}
 = \overline{\bar\nu(\hat\Phi(\bar u))} = \nu(\hat\Phi(\bar u)^*).
$$
{\nobreak
It follows } that $\hat\Phi(u\circ \hat S) = \hat\Phi(\bar u)^*$, and similarly $\hat\Phi(u\circ \hat S) = \hat\Phi(\bar u)^*$. So we have in \eqref{DeltaPhi1} (using $\overline{u\otimes v}=\bar u\otimes \bar v$):
\begin{align*}
\hat \Delta (\Phi\mu)(v\circ \hat S\otimes u\circ \hat S)
&= \mu\big( \hat\Phi(\bar v)^* \cdot \hat\Phi(\bar u)^*\big)
= \mu\big( \big(\hat\Phi(\bar u) \cdot \hat\Phi(\bar v)\big)^*\big)
= \overline{\bar\mu\big( \hat\Phi(\bar u) \cdot \hat\Phi(\bar v)\big)}
\\&= \overline{(\Phi\times\Phi)(\bar\mu)\big( \bar u\otimes\bar v\big)}
= \overline{\hat\Delta(\Phi\bar\mu)\big( \bar u\otimes\bar v\big)}
= \hat\Delta(\Phi\bar\mu)^* ( u\otimes v)
\\& = \hat\Delta(\Phi(\bar\mu)^*) ( u\otimes v)
 = \hat \Delta \hat S(\Phi\mu)(u\otimes v).
\end{align*}
\epr

We have obtained, in particular, that $S\hat\Phi(\omega) = \hat\Phi(\bar\omega)^*$ for all $\omega\in \Mhat_{*}$.

We have proved that $\Mhat$ is a quantum semigroup with involution. With respect to the involution on $\Mhat_*$ defined as in Definition \ref{def-involution}, that is $\omega^*(\Phi(\mu))=\bar\omega\circ\hat S(\Phi(\mu)$ for $\omega\in\Mhat_*$, $\mu\in M_*$, the map $\hat\Phi$ is involutive, as \eqref{omega-involution} shows.

\bprop\label{hatPhi**}
For $\omega\in\Mhat_*$, we have $\omega\in \Mhat_{**}$ if and only if $\hat\Phi(\omega)^* \in \hat\Phi(\Mhat)$. In this case, $\omega^* = \hat\Phi^{-1}\big( \hat\Phi(\omega)^*\big)$.
\eprop
\bpr
Take $\omega\in \hat M_{**}$, $\mu\in M_{**}$. Then:
$$
\bar\omega(\hat S\Phi\mu) = \overline{\omega( (\hat S\Phi\mu)^*)}
 = \overline{\omega( \Phi(\bar\mu))} = \overline{V(\omega,\bar\mu)}
 = \overline{\hat V(\bar\mu,\omega)} = \overline{\hat\Phi(\omega)(\bar\mu)}
 = \hat\Phi(\omega)^*(\mu).
$$
If $\omega^*$ is well defined, then the value above is also equal to $\omega^*(\Phi\mu) = \hat\Phi(\omega^*)(\mu)$, what proves the proposition.
\epr

The reasoning above is summarized in

\bprop\label{Mhat-is-QSI-dense}
Let $M$ be a quantum semigroup with involution, such that $M_{**} $ is dense in $M_*$. Then $\Mhat$ is a quantum semigroup with involution.
\eprop

\subsection{The case of $M_{**}$ not dense in $M_*$}

For general $M$, set $\Mhat=\hat{M_r}$. Proposition \ref{Mhat-is-QSI-dense} implies imediately

\btm\label{Mhat-is-QSI}
Let $M$ be a quantum semigroup with involution. Then $\Mhat$ is a quantum semigroup with involution.
\etm

Note however that $\Phi_M$ might not be a homomorphism if $M_{**}$ is not dense in $M_*$.

\bprop\label{hat-Phi-unitary}
If $\Mhat\ne\{0\}$ then $\hat\Phi: \Mhat_*\to ZM$ is unitary.
\eprop
\bpr
Since $Z:M_r\to M$ is a *-homomorphism, $\hat\Phi_M=Z_M\hat\Phi_{M_r}$ is a *-homomorphism even if $M\ne M_r$. If $V_r\in \Mhat\otimes M_r$ is the generator of $\Phi_{M_r}$ (which is unitary), then $\hat V = (Z_M\otimes 1)\big( \theta(V_r) \big) \in (Z_M M)\otimes\Mhat$, with the flip $\theta$, is the generator of $\hat\Phi_M$.
\epr

\bprop
$\hat\Phi$ is injective. 
\eprop
\bpr
If $\hat\Phi(\omega_1)=\hat\Phi(\omega_2)$ for some $\omega_1$, $\omega_2\in\Mhat_*$, then $\omega_1(\Phi(\mu))=\omega_2(\Phi(\mu))$ for every $\mu\in M_*$, so that $\omega_1=\omega_2$ on $\Phi(M_*)$. By weak density of $\Phi(M_*)$ in $\Mhat$, it follows that $\omega_1=\omega_2$.
\epr

\subsection{The axiom of the antipode}\label{sec-antipode}

In the theory of Hopf algebras, the antipode $S$ of a Hopf algebra $\cal M$ satisfies the axiom
\begin{equation}\label{ax-antipode}
 {\mathfrak m }(\Id\otimes S)\Delta = {\mathfrak m}(S\otimes \Id) \Delta = \varepsilon 1,
\end{equation}
where $\mathfrak m : \cal M\otimes\cal M\to \cal M $ is the multiplication and $\varepsilon: \cal M\to\C $ the counit.

Suppose that $\{0\}\ne M = \hat N $ and set $\cal M = \Phi_N(N_*) $. On $\cal M$, a counit is well defined (and satisfies the usual axioms of a counit): $\e\big(\Phi(\nu)\big) = \nu(1_N)$, $\nu\in N_*$.
Suppose that $N\subset B(H)$ is in its standard form, then every $\nu\in N_*$ has form $\nu=\nu_{xy}: a\mapsto \langle ax,y\rangle$, $a\in N$, $x,y\in H$. In this case, $\e(\nu_{xy}) = \langle x,y\rangle$.
For every $\nu\in N_*$, $\Phi(\nu)$ is a coefficient of the unitary representation $\hat\Phi_N: M_*\to N$. If we fix a basis $(e_\a) $ in $ H$, then in particular $\Phi(\nu_{e_\a,e_\beta}) = \hat\Phi_{\beta\a}$ in the notations \ref{nots-coeffs} for every $\a,\beta$.

Since $\Phi_N$ is a unitary representation of $N_*$, the formula \eqref{def-standard} implies (by decomposing $x$ and $y$ in the basis) that
\beq\label{Phi-counit-eq1}
\sum_\gamma \Phi(\nu_{y,e_\gamma})^* \Phi(\nu_{x,e_\gamma}) = 
\langle x,y\rangle \,1_M = \e(\nu_{xy}) \,1_M.
\eeq
At the same time,
$$
\Delta\big( \Phi(\nu_{xy})\big) = \sum_\gamma \Phi(\nu_{e_\gamma,y})\otimes \Phi(\nu_{x,e_\gamma}).
$$
by Proposition \ref{Delta-on-coeffs}. Together with $S_M \Phi(\nu_{e_\gamma,y}) = \Phi(\nu_{y,e_\gamma})^*$ which is valid by \eqref{hatSPhi},
the equality \eqref{Phi-counit-eq1} takes form
$$
{\mathfrak m}(S\otimes \Id)\Delta (a) =  \e(a) \, 1_M
$$
for $a=\nu_{xy}$. Similarly, on arrives at the other equality in \eqref{ax-antipode}. We see that the axiom of the antipode is satisfied for all $a\in\cal M$, in the sense described above. This suggests that every dual $M=\hat N$ is a quantum {\it group\/} and not just a semigroup.

\section{Morphisms and second duals}

In this section we do not suppose that $M_{**}$ is dense in $M_*$.

\begin{definition}
Let $\cal {QSI}$ be the category of quantum semigroup with involution. A morphism in $\QSI$ is a normal *-homomorphism $\phi: M\to N$ such that: $\Delta_N\circ\phi=(\phi\otimes\phi)\Delta_M$, $\phi(D(S_M))\subset D(S_N)$ and $S_N\circ\phi=\phi\circ S_M$. Note that we do not require that $\phi(1)=1$.
\end{definition}

\bprop\label{phi_*-involutive}
If $\phi: M\to N$ is a morphism in $\QSI$, then $\phi_*: N_*\to M_*$ is a *-homomorphism. In particular, $\phi_*(N_{**})\subset M_{**}$.
\eprop
\bpr
By definition, $\phi$ is ultraweakly continuous, so it has a pre-adjoint $\phi_*:N_*\to M_*$. Since $\phi$ is a coalgebra morphism, it is standard to show that $\phi_*$ is a homomorphism. It remains to prove that $\phi_*$ is involutive.

For any $\nu\in N_*$, $\overline{ \phi_*(\nu)}=\phi_*(\bar\nu)$, since for any $x\in M$
$$
\overline{ \phi_*(\nu)}(x) = \overline{ \phi_*(\nu)(x^*)} = \overline{ \nu\big(\phi(x^*)\big)}
 = \overline{ \nu\big(\phi(x)^*\big)} = \bar\nu\big(\phi(x)\big) = \phi_*(\bar\nu)(x).
$$
If moreover $\nu\in N_{**}$, then for any $x\in D(S_M)$
\begin{align*}
\overline{ \phi_*(\nu)}\circ S_M(x) &= \phi_*(\bar\nu)\big( S_M(x)\big)
 = \bar\nu\big( \phi( S_M(x))\big)
 = \bar\nu\big( S_N \phi( x)\big)
 = \nu^* \big( \phi( x)\big)  = \phi_* (\nu^*) (x).
\end{align*}
This calculation shows that $\phi_*(\nu)^*$ is well defined and equals $\phi_* (\nu^*)$.
\epr

Below we use the notations of Proposition \ref{quotient-by-perp}.

\bprop
If $\phi: M\to N$ is a morphism in $\QSI$, then there exists a morphism $\phi_r: M_r\to N_r$ such that $\phi_r Q_M = Q_N \phi$.
\eprop
\bpr
Recall that $\ker Q_M=(M_{**})^\perp$, $\ker Q_N=(N_{**})^\perp$. By Proposition \ref{phi_*-involutive}, $\phi_*(N_{**})\subset M_{**}$. It follows immediately that $\phi(\ker Q_M)\subset\ker Q_N$, so that $\phi_r$ which satisfies the equality $\phi_r Q_M = Q_N \phi$ is well defined. It is immediate to check that it is normal.

The fact that $\phi_r$ is a coalgebra morphism follows from the calculation below:
\begin{align*}
\Delta_{N_r}\raise1pt\hbox{$\phi_r$}\, Q_M &= \Delta_{N_r}Q_N \,\raise1.5pt\hbox{$\phi$} = (Q_N\otimes Q_N) \Delta_N\, \raise1.5pt\hbox{$\phi$}
 = (Q_N\otimes Q_N) (\phi\otimes\phi) \Delta_M
 \\& = (\phi_r Q_M\otimes \phi_r Q_M) \Delta_M
 = (\phi_r\otimes\phi_r) \Delta_{M_r} Q_M.
\end{align*}
For the coinvolution domains we have:
$$
\phi_r(D(S_{M_r})) = \phi_r( Q_M(D(S_M))) = Q_N\phi(D(S_M)) \subset Q_N D(S_N) = D(S_{N_r}),
$$
and moreover
\begin{align*}
\phi_r S_{M_r} Q_M = \phi_r Q_M S_M = Q_N \phi S_M = Q_N S_N\phi = S_{N_r} Q_N \phi = S_{N_r} \phi_r Q_M.
\end{align*}
\epr

\bprop\label{dual-morphism}
Let $\phi: M\to N$ be a morphism in $\cal {QSI}$. Then there is a dual morphism $\widehat\phi:\hat N\to \hat M$ such that $\hat\phi\circ\Phi_N=\Phi_M\circ\phi_*$ on $N_{**}^0$. If $\hat N\ne\{0\}$, then the equality holds on $N_{**}$.
\\
Moreover, $\hat\phi_r = \hat\phi$.
\eprop
\bpr
The statement is trivial if $\hat N=\{0\}$ or $\Mhat=\{0\}$ so we can assume that $N_{**}^0\ne\{0\}$ and $\Mhat\ne \{0\}$.

Consider $\Phi_M\circ\phi_*: N_{**}^0\to \Mhat$. Since $\Phi_M, \phi_*$ are *-homomorphisms on $M_{**}$ and $N_{**}$ respectively and $\phi_*(N^0_{**})\subset M_{**}$, their composition is a *-homomorphism on $N^0_{**}$, so it is lifted to a normal *-homomorphism $\hat\phi: \hat N\to \hat M$ such that $\hat\phi\,\Phi_N=\Phi_M\,\phi_*$ (on $N_{**}^0$). Both maps are in fact defined on $N_*$; by Proposition \ref{M*-hom-to-vna}, the equality holds on $N_{**}$, and if $N_{**}$ is dense in $N_*$, on $N_*$.

Now we can prove that $\hat\phi = \hat\phi_r$. By density, it is sufficient to prove that $\hat\phi \Phi_{N_r}  = \hat\phi_r \Phi_{N_r}$ on $(N_r)_{**}$, what is done as follows (recall Notations \ref{def-Z} and Definition \ref{def-Mhat-not-dense}): 
\begin{align*}
\hat\phi\, \Phi_{N_r} &= \hat\phi \,\Phi_{N_r} Z_{N*} Q_{N*} = \hat\phi \,\Phi_N Q_{N*} = \Phi_M \phi_* Q_{N*}
\\&
 = \Phi_M Q_{M*} \phi_{r*} = \Phi_{M_r} Z_{M*} Q_{M*} \phi_{r*} = \Phi_{M_r} \phi_{r*} = \hat\phi_r \Phi_{N_r}.
\end{align*}

From now on we can suppose that $\phi=\phi_r$, that is, $M_{**}$ is dense in $M_*$ and $N_{**}$ is dense in $N_*$.

Next we prove that $\hat\phi$ is a coalgebra morphism: $\Delta_{\Mhat} \, \widehat\phi = (\widehat\phi\otimes\widehat\phi)\Delta_{\widehat N}$. Since $\Delta_{\Mhat}$ is ultraweakly continuous, this equality is enough to check on $\Phi_N(N_{**}^0)$. Moreover, to check an equality in $\Mhat\bar\otimes\Mhat$ where $\hat\Delta_M$ takes its values, it is enough to consider evaluations on $x\otimes y$, with $x,y\in \Mhat_*$. We have, with any $\nu\in N_{**}^0$:
\begin{align}\label{nu-calculations}
\Delta_{\Mhat} \, \widehat\phi (\Phi_N(\nu))(x\otimes y) &= \widehat\phi(\Phi_N(\nu))(xy) = \Phi_M(\phi_*(\nu))(xy)
=\phi_*(\nu)(\hat\Phi_M(x)\hat\Phi_M(y))
\notag\\&=
\nu\Big(\phi(\hat\Phi_M (x))\,\phi(\hat\Phi_M (y))\Big).
\end{align}
By definition, $\Delta_{\widehat N}\Phi_N = \Phi_N\times\Phi_N$, so from the other side:
\begin{align*}
(\widehat\phi\otimes\widehat\phi)\Delta_{\widehat N}(\Phi_N(\nu))(x\otimes y)
&= (\Phi_N\times\Phi_N)(\nu)(\hat\phi_* x\otimes \hat\phi_*y)
= \nu\big(\hat\Phi_N(\hat\phi_* x)\cdot \hat\Phi_N(\hat\phi_*y)\big).
\end{align*}
Recalling that $\phi\hat\Phi_M = \hat\Phi_N\hat\phi_*$, we arrive at the required equality.

It remains to check the equality $S_{\widehat M}\widehat\phi = \widehat\phi S_{\widehat N}$. By definition $D(S_{\hat N}) = \Phi_N(N_*)$; for $\xi=\Phi_N(\nu)$, $\nu \in N_*$, we have $\hat\phi(\xi) = \hat\phi(\Phi_N(\nu)) = \Phi_M(\phi_*(\nu))\in \Phi_M(M_*) = D(S_{\Mhat})$. Then
$$
S_{\widehat M}\widehat\phi(\xi) 
 = S_{\widehat M}\Phi_M(\phi_*(\nu)),
$$
and by Definition \ref{def-hat-S}, for any $\omega\in \Mhat$:
\begin{align*}
\omega\big( S_{\widehat M}\Phi_M(\phi_*(\nu))\big) &= V_M^*(\omega,\phi_*(\nu)) = \overline{ V_M( \bar\omega, \overline{\phi_*(\nu)} )}
= \overline{ V_M( \bar\omega, \phi_*(\bar\nu) )}
\\& = \overline{\bar\omega\big( \Phi_M(\phi_*(\bar\nu))\big)}
= \overline{\bar\omega\big( \hat\phi (\Phi_N(\bar\nu))\big)}
= \overline{\hat\phi_*(\bar\omega)\big( \Phi_N(\bar\nu)\big)}
\\&= \overline{V_N\big( \hat\phi_*(\bar\omega), \bar\nu\big)}
= \overline{V_N\big( \overline{\hat\phi_*(\omega)}, \bar\nu \big)}
= V_N^* \big( \hat\phi_*(\omega), \nu \big)
\\&= \hat\phi_*(\omega)\big( S_{\hat N}\Phi_N( \nu ) \big)
= \omega\big( \hat\phi\,S_{\hat N}(\Phi_N( \nu )) \big)
= \omega\big( \hat\phi\,S_{\hat N}(\xi) \big).
\end{align*}
It follows that $S_{\widehat M}\widehat\phi(\xi) = \hat\phi\,S_{\hat N}(\xi)$ as required.
\epr

\brem\label{hat-is-functor}
If we have two $\cal {QSI}$-morphisms $\phi:M\to N$ and $\psi: N\to L$, then on $L_{**}^0$ we have:
$$
\hat{\psi\circ\phi}\circ\Phi_L = \Phi_M\circ (\psi\circ\phi)_* = \Phi_M\circ\phi_*\circ\psi_*,
$$
$$
\hat\phi\circ\hat\psi \circ\Phi_L = \hat\phi \circ \Phi_N \circ \psi_*.
$$
It might happen that $\psi_*(L^0_{**})$ is not a subset of $N^0_{**}$ so we cannot continue the last line as $\Phi_M\circ\phi_*\circ\psi_*$, in general.
If $\hat N\ne\{0\}$, then we can use the fact that $\hat\phi \circ \Phi_N=\Phi_M\circ\phi_*$ on $N_{**}$ and conclude that both displayed lines are equal, what implies $\hat{\psi\circ\phi}=\hat\phi\circ\hat\psi$. If $\hat N=\{0\}$, then necessarily $\hat\phi\circ\hat\psi=0$ but it might happen that $\hat{\psi\circ\phi}\ne0$.
\erem

\bprop
Let $\Id_M$ be the identity morphism of $M$, then $\hat\Id_M=\Id_{\Mhat}$.
\eprop
\bpr
One checks trivially that $(\Id_M)_*=\Id_{M*}$. 
Next, by definition $\hat\Id_M \Phi_M = \Phi_M (\Id_M)_*$, and then $\hat\Id_M \Phi_M = \Phi_M$ on $M^0_{**}$. Since $\Phi_M(M^0_{**})$ is dense in $\hat M$, it follows that $\hat\Id_M = \Id_{\Mhat}$.
\epr

\bprop\label{D-exists}
For every $M$, there is a normal *-homomorphism $D_M: \widehat{\Mhat}\to M$ such that $D_M\circ\Phi_{\Mhat}(x)=\hat\Phi_M(x)$ for all $x\in \Mhat^0_{**}$. If $\hat{\Mhat}\ne\{0\}$, then this equality holds actually on $\Mhat_{**}$. If $M_{**}$ is dense in $M_*$, then $D_M$ is a coalgebra morphism and thus a morphism in $\cal {QSI}$.
\eprop
\bpr
If $\Mhat=\{0\}$, the statement holds for $D_M=0$ which is a morphism; below we suppose that $\Mhat\ne\{0\}$.

Denote $N=\Mhat$. The canonical map $\hat \Phi_M: N_*\to M$ is a *-homomorphism. In particular, $\hat\Phi_M|_{N_{**}^0}$ is a *-homomorphism, so it is extended uniquely to a normal homomorphism of its von Neumann envelope: $D_M:\widehat N\to M$. By definition, $D_M$ satisfies the equality in the statement for $x\in N_{**}^0=\Mhat^0_{**}$. If $\hat{\Mhat}\ne\{0\}$, then by Propositions \ref{hat-Phi-unitary}, \ref{unitary-is-non-degenerate} and \ref{M*-hom-to-vna} this equality holds also on $\Mhat_{**}$.

Similarly, there exists a unique $D_{M_r}: \hat{\Mhat}\to M_r$ such that $D_{M_r}\circ\Phi_{\Mhat} = \hat\Phi_{M_r}$ on $\Mhat_{**}$. Recalling that $\hat\Phi_M = Z_M \hat\Phi_{M_r}$, we infer that $D_M=Z_M D_{M_r}$.

To prove the last statement, we can suppose that $M_{**}$ is dense in $M_*$ and $Z_M=1$. In this case $D_M\circ\Phi_{\Mhat}=\hat\Phi_M$ on $\Mhat_*$, and the preadjoint map $(\hat\Phi_M)_*: M_*\to \hat M$, which is by definition $\Phi_M$, equals also to $(\Phi_{\Mhat})_*\circ(D_M)_* = \hat\Phi_{\Mhat}\circ (D_M)_*$. Also, $\Phi_M$ is a homomorphism.

To prove that $D_M$ is a morphism of coinvolutive Hopf--von Neumann algebras, one should check the equality
\begin{align*}
\Delta_M \, D_M (x)(\mu\otimes \nu) = (D_M\otimes D_M)\Delta_{\widehat N}(x)(\mu\otimes \nu)
\end{align*}
for every $x\in \hat N$ and $\mu,\nu\in M_*$. By density, it is sufficient to consider $x=\Phi_N(y)$ with $y\in N_{**}^0$. Then we have:
\begin{align*}
\Delta_M \, D_M\Phi_N(y)(\mu\otimes \nu) &= \Delta_M \, \hat\Phi_M(y)(\mu\otimes \nu)=
\hat\Phi_M(y)(\mu\nu) = y(\Phi_M(\mu\nu))
\\&= y( \Phi_M(\mu)\Phi_M(\nu)) = y(\hat\Phi_N((D_M)_*(\mu))\cdot\hat\Phi_N((D_M)_*(\nu)))
\\&=\Phi_N(y)((D_M)_*(\mu)\cdot(D_M)_*(\nu)) = \Delta_{\widehat N}(x)((D_M)_*(\mu)\otimes(D_M)_*(\nu))
\\&= (D_M\otimes D_M)(\Delta_{\widehat N}(x))(\mu\otimes \nu).
\end{align*}

To verify that $D_M$ agrees with the coinvolutions, let $y\in \Mhat_*$. Then
$D_M\Phi_{\Mhat}(y) = \hat\Phi_M(y) \in D(S_M)$ by Proposition \ref{hatPhi(Mhat)-in-D(S)}, and $S_M\big( \hat\Phi_M(y) \big) = \hat\Phi_M(\bar y)^*$ by formula \eqref{hatSPhi}. From the other side,
\begin{align*}
D_M S_{\hat {\Mhat}} \Phi_{\Mhat}(y)&= D_M \big( \Phi_{\Mhat}(\bar y)^* \big)
 = \big( D_M\Phi_{\Mhat}(\bar y)\big)^*
 = \big( \hat\Phi_M(\bar y)\big)^*
 = S_M\big( \hat\Phi_M(y) \big)
 = S_M D_M\Phi_{\Mhat}(y).
\end{align*}
This proves that $S_M D_M= D_MS_{\widehat {\Mhat}}$ on $\Phi_{\Mhat}(\Mhat_*)$.
If $S_{\hat{\Mhat}}$ is extended by Proposition \ref{S-extended}, then by continuity (Proposition \ref{R-continuous}) we still have equality on the whole of $D(S_{\hat{\Mhat}})$.
\epr

As it was shown already in the case of a bounded $S$ in \cite{haar}, the map $D_M$ need not be neither injective nor surjective. For the algebra $M$ from the example 5.10 of \cite{haar}, $\hat\Mhat=0$, so $D_M=0$. For a second example, take $M=L_\infty(G)$. Then $\hat\Mhat=C_0(G)^{**}$, and $D_M$ is a quotient map but is not injective.

If $M$ is a dual of another algebra, then $D_M$ is right invertible:

\bprop\label{DE=id}
If $M=\hat N$ for some $N$ and $\hat{\Mhat}\ne\{0\}$, then there is a $\cal {QSI}$-morphism $E_M: M\to \widehat {\widehat M}$, such that $D_M\circ E_M=\id_M$.
\eprop
\bpr
Since $M=\hat N = \hat {N_r}$, we can assume that $N=N_r$.
By Proposition \ref{D-exists} applied to $N$, since $\Mhat\ne\{0\}$, there is a $\cal {QSI}$-morphism $D_N: \widehat{\hat N}\to N$ such that $D_N\circ\Phi_{\hat N}(x)=\hat\Phi_N(x)$ for all $x\in \hat N_{**}=M_{**}$. By Proposition \ref{dual-morphism}, there exists a dual morphism $E_M=\widehat D_N: \hat N=M\to (\widehat{\hat N})\widehat{\;} = \widehat{\widehat M}$.

By assumption $\Mhat=\hat{\hat N}\ne0$, and in this case it was proved in Proposition \ref{D-exists} that $\Phi_N = \hat\Phi_M\circ (D_N)_*$ on $N_*$. Moreover, since $\hat{\Mhat}\ne\{0\}$, we have $D_M\circ \Phi_{\hat M}=\hat\Phi_M$ on $\Mhat_{**}$.

By Proposition \ref{dual-morphism}, since $M\ne\{0\}$, we have $E_M\circ\Phi_N=\Phi_{\hat M}\circ (D_N)_*$ on $N_{**}$. Then on $N_{**}$,
$$
D_M\circ(E_M\circ\Phi_N) = (D_M\circ \Phi_{\hat M})\circ (D_N)_* = \hat\Phi_M\circ (D_N)_* = \Phi_N;
$$
note that $(D_N)_*(N_{**})\subset \Mhat_{**}$ by Proposition \ref{phi_*-involutive}.
Thus, $D_M\circ E_M=\id_M$ on $\Phi_N(N_{**})$. Since this latter is weakly dense in $M$, this equality holds everywhere.
\epr

\bcor\label{M**-is-dense}
If $M=\hat N$ for some $N$ and $\hat{\Mhat}\ne\{0\}$, then $M_{**}$ is dense in $M_*$.
\ecor
\bpr
In Proposition \ref{D-exists} it was proved that $D_M=Z_M D_{M_r}$. By Proposition \ref{DE=id}, $D_M$ is surjective. This implies that $\zeta_M=1_M$ and $M_{**}$ is dense in $M_*$.
\epr

\bprop\label{E-unital}
If $M=\hat N$ for some $N$ and $\Mhat^{(4)}
\ne\{0\}$ then $E_M$ is unital.
\eprop
\bpr
It is clear that $E_M(1)=p=p^2=p^*$. Suppose that 
$p\ne1$. Then $E_M(M)$ is contained in the proper weakly closed *-subalgebra $I:=p\hat{\hat M}p$.
Consider $\pi=\Phi_M\circ E_{M*}:\hat{\hat M}_*\to \hat M$. Since $M_{**}$ is dense in $M_*$, $\Phi_M$ is a *-homomorphism; by Proposition \ref{phi_*-involutive}, $\pi$ is a *-representation. For $\omega\in \hat M_*$ , $\omega(\pi(\mu)) = \Phi_M\circ E_{M*}(\mu)(\omega) = \mu(E_M\circ\hat\Phi_M(\omega))$, so $\pi_*(\omega)=E_M\circ\hat\Phi_M(\omega)$ and the space of coefficients of $\pi$ is contained in the subalgebra $I$. Then the equality \eqref{def-standard} cannot hold, so $\pi$ is not unitary. By Corollary \ref{M**-is-dense}, $\hat{\Mhat}_{**}$ is dense in $\hat{\Mhat}_*$, so we can apply Proposition \ref{nondegenerate-is-standard}; it follows that $\pi$ is degenerate on $\hat{\hat M}{}_*^0$, or equivalently $\pi(\hat{\hat M}{}_*)$ is not contained in the weak closure of $\pi(\hat{\hat M}{}_{**}^0)$. In particular, $\pi(\hat{\hat M}{}_{**}^0)$ is not weakly dense in $\Mhat$.

Consider now $\hat E_M:\hat{\hat {\hat M}}\to \hat M$. Since $\hat{\hat {\hat M}}\ne\{0\}$, by Proposition \ref{dual-morphism} $\hat E_M\Phi_{\hat{\hat M}}=\Phi_ME_{M*}=\pi$ on $\hat{\hat M}{}_{**}$. and since $\hat E_M$ is weakly continuous, $\hat E_M(\hat{\hat {\hat M}})$ is contained in the closure of $\hat E_M(\Phi_{\hat{\hat M}}(\hat{\hat M}{}_{**}^0))= \pi(\hat{\hat M}{}_{**}^0)$, so $\hat E_M$ is not surjective. From the other side, since $\hat{\hat {\hat M}}\ne\{0\}$, we have by Remark \ref{hat-is-functor} equality $\hat E_M\hat D_M=(D_ME_M)\hat{\;}={\hat{\Id}}_M=\Id_{\Mhat}$ so $\hat E_M$ must be surjective. This contradiction proves that $E_M$ is in fact unital.
\epr

\bprop
$(D_M)_*(M_*)$ is a two-sided module over $\hat{\Mhat}$.
\eprop
\bpr
We write further $D$, $E$ for $D_M$, $E_M$ respectively.
First, $D_*(M_*) = (\ker D)_\perp$ in $\hat{\Mhat}_*$: For $\omega\in \hat{\Mhat}_*$ if $\omega=D_*(\mu)\in D_*(M_*)$ then $\omega(x) = \mu(Dx)=0$ for $x\in \ker D$; if $\omega\in(\ker D)_\perp$ then set $\mu=E_*(\omega)$; for $x\in \hat{\Mhat}$, as $x-EDx\in \ker D$, we have $\omega(x) = \omega(EDx) = \mu(Dx)$, so that $\omega=D_*(\mu)\in D_*(M_*)$.

Now $(x.\omega)(y)=\omega(yx)=0$ and $(\omega.x)(y)=\omega(xy)=0$ for all $\omega\in D_*(M_*)$, $x\in \hat{\Mhat}$, $y\in \ker D$, what proves the proposition.
\epr

\section{Locally compact quantum groups
}\label{sec-LCQG}

Every von Neumann algebraic quantum group is a quantum semigroup with involution in our definition (see for example \cite{kust-vaes-vna}). Let us compare our construction with the universal dual of Kustermans \cite{kustermans}.

The Kustermans' universal quantum group $A_u$ is defined for a $C^*$-algebraic locally compact quantum group $A$. If $M$ is a von Neumann algebraic locally compact quantum group, then \cite{kust-vaes-vna} it contains a canonically defined $C^*$-subalgebra $A$ which is a $C^*$-algebraic locally compact quantum group with inherited structure.

The space $L^1(A)$ used as a starting point in \cite{kustermans} is isomorphic to $M_*$ \cite[p.913]{kust-vaes}. The subspace $L^1_*(A)$ which carries the involution is exactly our $M_{**}$, and is dense in $M_*$ \cite[p.303]{kustermans}. By \cite[Corollary 4.3]{kustermans}, for every non-degenerate *-representation $\pi: L^1_*(A)\to B(H)$ there exists a unitary $U\in B(H)\otimes M(A)$ such that $\pi(\mu)=(\Id\otimes\mu)(U)$ for all $\mu\in L^1_*(A)$. Since $M(A)\subset A''=M$, 
it follows that $U\in B(H)\otimes M$ and $\pi$ extends to a unitary representation of $M_*$.

This implies that $M_{**}=M^\times_{*}=M^0_{**} = L^1_*(A)$, and $\hat A_u=C^*(L^1_*(A))=C^*(M_{**}^0)$. Our dual algebra $\Mhat=\hat A_u^{**}$ is the enveloping von Neumann algebra of the universal dual $\hat A_u$.

\btm\label{th-LCQG}
Let $A$ be a $C^*$-algebraic locally compact quantum group and $M=A''$ its associated von Neumann algebraic locally compact quantum group. Then:
\begin{enumerate}
\item $\Mhat = \hat A_u^{**}$;
\item $\hat{\Mhat} = A_u^{**}$.
\end{enumerate}
In particular, $\hat{\hat{\Mhat}} = \Mhat$.
\etm
\bpr
(1) follows from the discussion above.

(2)
Since $\hat{\Mhat}=C^*(\Mhat^0_{**})^{**}$ and $A_u = C^*(L_1^*(\hat A))$, the question is to show that $C^*(\Mhat^0_{**}) = C^*(L_1^*(\hat A))$. The epimorphism $\hat\pi: \hat A_u\to\hat A$ \cite[2.15]{kustermans} generates an epimorphism $\hat\pi^{**}: \hat A_u^{**}\to \hat A^{**}$ and an imbedding $\hat\pi^*: \hat A^*\to \hat A_u^{*} = \Mhat_*$.

The equality $(\imath\otimes\hat\Delta_u)(\hat{\cal V}) = \hat{\cal V}_{13}\hat{\cal V}_{12}$ \cite[3.5]{kustermans} implies that $\lambda_u^*$ is anti-multiplicative.
Let us prove now that $\hat\pi^*(\hat A^*)$ is an ideal in $\hat A_u ^*$.
It is sufficient to show that for any $\alpha\in \hat A^*$ and any representation $\rho: \hat A_u\to B(H)$ and every $\omega\in B(H)_*$ $\hat\pi^*(\alpha)\rho^*(\omega)\in \hat\pi^*(\hat A^*)$. For $\mu\in L^*_1(A)$:
$$
\big(\hat\pi^*(\alpha)\rho^*(\omega)\big) \big(\lambda_u(a)\big)
 = \lambda_u^*\big(\hat\pi^*(\alpha)\rho^*(\omega)\big) (a)
 = \big( \lambda_u^*\rho^*(\omega)\lambda_u^*\hat\pi^*(\alpha)\big) (a)
 = \big( \lambda_u^*\rho^*(\omega)\lambda^*(\alpha)\big) (a)
 = \lambda^*(\xi) (a),
$$
since $\lambda^*(\hat A^*)$ is an ideal in $A$.

Every non-degenerate $*$-representation of $L_1^*(\hat A)$ on a Hilbert space $H$ is extended to a unitary representation of $L_1(\hat A)$ with a generator $U\in B(H)\otimes M(\hat A)$ \cite[Corollary 4.3]{kustermans}. By \cite[Proposition 3.13]{kustermans}, there is unique $V_u\in B(H)\otimes M(\hat A_u)$ such that $(\hat\pi^{**}\otimes\Id)(V_u)=U$ and the corresponding map of $L_1(\hat A_u)$ is a *-representation. Inclusion $M(\hat A_u)\subset \hat A_u^{**}=\Mhat$ implies that $V\in B(H)\otimes \Mhat$, and $(\Id\otimes \hat\pi^{**})(V_u) = U$. Thus, every non-degenerate $*$-representation of $\hat\pi^*\big(L_1^*(\hat A)\big)$ is extended to a unitary representation of $\Mhat_*$, and since such the extension is unique, it follows that every non-unitary representation vanishes on $L_1^*(\hat A)$ and $L_1^*(\hat A)\subset \Mhat^\times_{**}$.

From the other side, if $\rho: \Mhat_*\to B(H)$ vanishes on $L_1^*(\hat A)$ then its coefficients are contained in $L_1^*(\hat A)^\perp = L_1(\hat A)^\perp \subset \Mhat$.
Let us show that $L_1(\hat A)^\perp$ (recall that we identify $L_1(\hat A)$ with $\hat\pi^*(L_1(\hat A))$) is a proper weakly closed ideal, then it will follow that $\rho$ is non-unitary.

The epimorphism $\hat\pi: \hat A_u\to \hat A$ is extended, by the universality property, to a $*$-homomorphism $\tilde\pi: \Mhat=\hat A_u^{**}\to \hat A''$. Then, since $L_1(\hat A)=(\hat A'')_*$, one has  
$\Big(\tilde\pi^*(L_1(\hat A))\Big)^\perp = \ker\tilde\pi$. Since $\mu(\tilde\pi(x))=\mu(\hat\pi(x))$ for $x\in \hat A_u$ and $\mu\in L^1(\hat A)$, by weak density of $\hat A_u$ in $\Mhat$ it follows that $\tilde\pi^*(\mu)=\hat\pi^*(\mu)$, so that $\hat\pi^*(L_1(\hat A))^\perp= \ker\tilde\pi$.
This proves (2).
\epr

Theorems 7.4 and 7.6 of \cite{haar} are valid for our case also. The proofs are identical, with replacement of $M^0_*$ by $M^0_{**}$, and adding the remark of Proposition \ref{involutive=coeff_in_D(S)}. We will not repeat the proofs and just state the results:

\btm\label{commutative-is-group}
Let $M\ne\{0\}$ be commutative and $M\simeq \hat N$ for some $N$. Then there is a locally compact group $G$ such that $M\simeq C_0(G)^{**}$.
\etm

\bcor
If $\{0\}\ne M\simeq \hat{\hat M}$ and $M$ is commutative, then $M\simeq C_0(G)^{**}$ for a locally compact group $G$.
\ecor

\btm\label{cocommutative-is-group}
Let $M$ be cocommutative (i.e. $M_*$ is commutative) and $M\simeq \hat N$ for some $N$. Then, if $\Mhat\ne\{0\}$, there is a locally compact group $G$ such that $M\simeq W^*(G)$.
\etm

\section{examples}\label{sec-examples}

Outside the class of locally compact quantum groups, some examples of algebras and their duals are given in \cite[Examples 5.10--5.12, 8.8]{haar}. The present construction extends the one of \cite{haar}, so the same examples are valid for it. They give several algebras $M$, commutative or not, such that $\Mhat$ is $\{0\}$ or $\C$. In \cite[Example 5.11]{haar}, the structure of $L^\infty(\R^2) $ is changed in such a way that its dual becomes isomorphic to $L^\infty(\R) $. This demonstrates the idea that our duality cuts out the ``non-unitary'' part of a given algebra, and leaves the ``unitary'' one.

Below are presented some more examples. In every one constructed up to now, the dual $\Mhat$ coincides with the dual of some locally compact quantum group.

\begin{example}
Let $B$ be the quantum semigroup $C({\widetilde S}^+_N)$ defined by Banica and Skalski \cite{banica-skalski}. Recall that it is defined starting with a ``submagic'' $N\times N$ matrix $u=(u_{ij})$ with entries in a unital $C^*$-algebra $A$. Being `submagic'' means that $u_{ij}=u_{ij}^*=u_{ij}^2$ for every $i,j$, and $u_{ij}u_{ik}=u_{ji}u_{ki}=0$ for every $i$ if $j\ne k$. By definition, $B$ is the universal unital $C^*$-algebra with the relations above. The authors show that $B$ admits a comultiplication (a unital coassociative *-homomorphism) defined by the formula $\Delta(u_{ij}) = \sum_k u_{ik}\otimes u_{kj}$, and a ``sub-coinvolution'' (everywhere defined *-antihomomorphism) defined by $S(u_{ij})=u_{ji}$.

Let $M$ be the enveloping von Neumann algebra of $B$, then $\Delta$ and $S$ extend obviously to normal maps on $M$. It is immediate to verify that $S$ is a proper coinvolution on $M$. Since $S$ is bounded, $M_{**}=M_*$.
The elements $(u_{ij})$ are coefficients of a *-representation of $M_*$ (by Theorem \ref{rep-by-coeffs} for example). One shows easily that it is irreducible.

However, $u$ is clearly non-unitary, thus every $u_{ij}$ belongs to the annihilator of $M^\times_{**}$ (see Definition \ref{def-Mtimes}). Since $(u_{ij})$ generate $B$, by Definition \ref{def-M0**} $M^0_{**}=\{0\}$ and $\Mhat=\{0\}$.

\end{example}

\begin{example}
X. Li \cite{Li} defines a reduced $C^*$-algebra $C_r^*(P)$ of a discrete left cancellative semigroup $P$, using its regular representation on $\ell^2(P)$. The algebra $C_r^*(P)$ is generated by the translation operators $T_p$, $p\in P $, and their adjoints $T_p^*\in B(\ell^2(P))$.
As a linear space, $C_r^*(P)$ is generated by $E_X L_g$, where $X$ is an ideal in $P$, $E_X$ is the operator of multiplication by the characteristic function of $X$, and $g=p_1^{\pm1}\dots p_n^{\pm1}$ with $p_j\in P$.

Set $M = C^*_r(P)^{**}$. Li does not define a coinvolution on his algebra. In order that it fits into our assumptions, set $S(T_p)=T_p^*$, and accordingly $S(T_p^*) = T_p$, $p\in P$, and extend it as a linear anti-homomorphism onto the algebra generated by these elements. For every $f\in \ell^2(P)$ and every $g_j\in G$, $X_j\subset P$, $j=1,\dots,n$ the function $\sum \bar\lambda_j E_{X_j} L_{g_j} f$ is the complex conjugate of $\sum \lambda_j E_{X_j} L_{g_j} \bar f$, and as $\|f\|=\|\bar f\|$, this implies that $S$ is isometric (and bounded) on $C_r^*(P)$.

By definition $\Delta(T_p)=T_p\otimes T_p$, so that every $T_p$ is a character of $M_*$. With the coinvolution above, it is involutive. It is unitary as an element of $M$ if and only if $p$ is invertible in $P$.

By construction, $\Delta(E_X)=E_X\otimes E_X$ for an ideal $X$ in $P$, so that it is also a character, which is clearly unitary only if $X=\emptyset$ or $X=P$. It follows that $M^\times_{**}$ is the linear dual space of ${\rm lin}\{L_g: g\in P\cap P^{-1}\}$. As $H=P\cap P^{-1}$ is a group, $M^\times_{**}$ is isomorphic to $\ell^1(H)^*$. It is readily seen that in fact, $M^0_{**}=M^\times_{**}$.
We conclude that $\Mhat = C^*(H)^{**}$.
\end{example}

\subsection{Weakly almost periodic compactifications}\label{dual-cstar}

\def\G{\mathbb{G}}

Let $G$ be a locally compact group, and let $P$ be its weakly periodic compactification. It is known (see Example \ref{ex-semigroup}) that $P$ is a compact semitopological semigroup, and $M=C(P)^{**}$ is a quantum semigroup with involution.

In $P$, there exists the minimal ideal $J$ (so that $PJ=JP=J$ and $J$ does not contain any other ideals), which has the form $J=eP$ for a central idempotent $e$ and is isomorphic (and homeomorphic) to the Bohr compactification $bG$ of $G$ \cite[Theorem III.1.9]{ruppert}.

Let $H$ be a locally compact group and let $\phi_H: H\to P$ be a homomorphism which is also a homeomorphism of $H$ onto $\phi_H(H)$. We are interested in two cases: $H=G$, with the canonical imbedding into $P$, and $H=bG$.
In this setting, let $p_H:C(P)\to C_b(H)$, $f\mapsto f\circ\phi_H$, be the ``restriction'' map, and let $\tau_H:M(H)\to M(P)$ be the ``extension'' map: $\tau_H(\mu)(f) = \int p_H(f) d\mu$ for $\mu\in M(H)$, $f\in C(P)$. Its dual $\tau^*_H: C(P)^{**}\to C_0(H)^{**}$ is ultraweakly continuous and extends $p$, so it is a *-homomorphism.

\bprop
If $\pi: M(P)\to B(L)$ is a unitary representation, then $\pi\circ\tau_H$ is a unitary representation of $M(H)$.
\eprop
\bpr
Fix a basis $(e_\a)$ in $L$. For every $\a,\beta$ we have $\pi_{\a\beta}\circ\tau \in C_0(H)^{**}$, and, by continuity,
\begin{align*}
\sum_\gamma (\pi_{\gamma\a}\circ\tau)^*\cdot\pi_{\gamma\beta}\circ\tau = \tau^*_H \Big( \sum_\gamma \pi_{\gamma\a}^*\cdot\pi_{\gamma\beta} \Big)
 = \begin{cases} 1, &\a=\beta\\ 0, &\a\ne\beta\end{cases}
\end{align*}
as well as $\sum_\gamma \pi_{\a\gamma}\circ\tau \cdot(\pi_{\beta\gamma}\circ\tau)^*$. This proves that $\pi\circ\tau_H$ is unitary.
\epr

\bprop\label{irreps-of-M(P)}
Every irreducible unitary representation $\pi$ of $M(P)$ is finite-dimensional.
\eprop
\begin{proof}
As usual, let $\delta_t$, $t\in P$, be the probability measure concentrated at $t$, and set $\bar\pi(t)=\pi(\delta_t)$. The equations \eqref{def-standard} imply that $\bar\pi(t)$ is a unitary operator for every $t\in P$. It follows that $\bar\pi: P\to B(L)$ is a representation, which is unitary when restricted
to $\phi_H(H)$, and the Proposition above implies \cite[Theorem 4.4]{haar} that $\bar\pi\circ\phi_H$ is (strongly) continuous, and $\pi\circ\tau_H(\mu) = \int_H \bar\pi\circ\phi_H(t) d\mu(t)$ for every $\mu\in M(H)$. It follows also that $\bar\pi\circ\phi_H$ and $\pi\circ\tau_H$ have the same invariant subspaces.

For the idempotent $e$, $\bar\pi(e)$ is unitary and idempotent, so it is the identity operator.
For every $\mu\in M(P)$, $\delta_e*\mu$ is in the image of the map $\tau_{bG}$ defined above. Indeed, $\delta_e*\mu = \tau_{bG}(\nu)$, where $\nu\in M(bG)$ is defined as $\nu(f) = \int_{P} f(\phi_{bG}^{-1}\big(e\phi_{bG}(t)\big) d\mu(t)$, $f\in C(bG)$. As $\pi(\mu) = \pi(\delta_e*\mu)$, we see that any $\tau_{bG}$-invariant subspace is also $\pi$-invariant.

Suppose now that $\pi$ is irreducible. Then, by discussion above, so is $\bar\pi\circ\phi_{bG}$, and being continuous, it is finite-dimensional.
\end{proof}

Conversely, if $\sigma$ is an irreducible continuous representation of $bG$, then $\pi(\mu) = \int_{P} \sigma(\phi_{bG}^{-1}(et)) d\mu(t)$ defines a representation of $M(P)$, and one easily verifies that the conditions \eqref{def-standard} hold, so that $\pi$ is unitary.

As shows the reasoning above, $I=\tau_{bG}(M(bG))$ is an ideal in $M(P)$. It follows that every irreducible representation $\pi$ either vanishes on $I$ or is irreducible on it. If $\pi$ is non-unitary and $\pi|_I$ is irreducibe, then $\pi\circ\tau_{bG}$ is non-unitary (otherwise $\pi$ would be unitary, as the unique extension from $I	$), so that $\pi\circ\tau_{bG}$ vanishes on $L^1(bG)$. We see that in any case $\pi$ vanishes on $\tau_{bG}(L^1(bG))$, thus the latter is contained in $M^\times_{**}$. The annihilator of $\tau_{bG}(L^1(bG))$ in $C(P)^{**}$ is an ultraweakly closed (proper) ideal, what can be proved in a virtually the same way as in \cite[Proposition 4.3]{haar}. This implies that $(M^\times_{**})^\perp\ne C(P)^{**}$ and as a consequence $M^0_{**}=M^\times_{**}$, according to Definition \ref{def-M0**}.

For a representation $\sigma$ of $M(bG)$, set $\pi_\sigma(\mu) = \sigma(\tau_{bG}^{-1}(\mu*\delta_e))$, $\mu\in M(P)$ (it is clear that $\tau_{bG}$ is injective). Proposition \ref{irreps-of-M(P)} above implies that for $\mu\in M^0_{**}$,
\begin{align*}
\sup\{ &\|\pi(\mu)\|: \pi\text{ is a *-representation of }M^0_{**}\}
\\
&= \sup\{ \|\pi_\sigma(\mu)\|: \sigma\text{ is a unitary *-representation of }M(bG)\},
\end{align*}
so that $C^*(M^0_{**})$ is isomorphic to $C^*(M(bG)^0)$. As we know, this equals to $C^*(bG)$, so finally  $\hat M = C^*(bG)^{**}$.

\end{document}